\documentstyle[]{amsart}
\numberwithin{equation}{section}

\newtheorem{theorem}[equation]{Theorem}
\newtheorem{lemma}[equation]{Lemma}
\newtheorem{corollary}[equation]{Corollary}
\newtheorem{proposition}[equation]{Proposition}
\theoremstyle{definition}
\newtheorem{definition}[equation]{Definition}
\newtheorem{example}[equation]{Example}

\theoremstyle{remark}

\newcommand{\Eff}{\operatorname{Eff}}

\newcommand{\Pic}{\operatorname{Pic}}

\newcommand{\CR}{\operatorname{CR}}
\newcommand{\qprod}{\operatorname{\ast_Q}}
\def\surf{X}
\def\surff{Y}

\newcommand{\import}{\input}

\begin{document}
\title{Quantum Cohomology of Rational Surfaces}
\author{Bruce Crauder}
\address{Dept. of Mathematics\\Oklahoma State
University\\Stillwater, OK 74078}
\email{crauder@@math.okstate.edu}
\author{Rick Miranda}
\address{Dept. of Mathematics\\Colorado State University\\Ft.
Collins, CO
80523}
\email{miranda@@riemann.math.colostate.edu}
\thanks{Research supported in part by the NSF under grant
DMS-9403392}
\thanks{\today. Set in type by AMSLaTeX Version 1.1}

\maketitle

\begin{abstract}
In this article formulas for the quantum product of a rational
surface are given, and used to give an algebro-geometric proof
of the associativity of the quantum product for strict Del Pezzo surfaces,
those for which $-K$ is very ample.  An argument for the associativity
in general is proposed, which also avoids resorting to the symplectic
category.  The enumerative predictions of Kontsevich and Manin
\cite{kontsevich-manin}
concerning the degree of the rational curve locus in a linear system
are recovered.  The associativity of the quantum product for the cubic
surface is shown to be essentially equivalent to the classical enumerative
facts concerning lines: there are $27$ of them, each meeting $10$ others.
\end{abstract}

\tableofcontents

\section*{Introduction}

The purpose of this article is to give an algebro-geometric
description
of the quantum cohomology ring for a general rational surface
$\surf$.
By a ``general'' rational surface we mean one in which all linear
systems
have the expected dimension, and in which the locus of rational
curves
in each linear system also has the expected dimension.
We understand that it is not known whether the general blowup of
the plane in ten or more general points
is general in this sense;
however we proceed anyway, developing the quantum cohomology ring
based on the expected linear systems.
In this sense the theory is entirely a numerical one.
We view the genericity assumption on the surface $\surf$
as capable of replacing the genericity assumptions
for the symplectic geometry on which the quantum theory is usually
based
(see for example \cite{ruan-tian,siebert-tian} and a forthcoming
paper
of Grassi \cite{grassi}
which also addresses aspects of the quantum cohomology of
rational surfaces).
It seems to us desirable to have a description of the quantum
product
which is based on algebraic geometry rather than symplectic
geometry,
and this is what we have tried to offer in the case of rational
surfaces.

In Section \ref{sectionX3class} we explain how classes in the
triple product $\surf^3$ are used to define classes on $\surf$
itself,
via the K\"unneth formula; this is formal.
In Section \ref{sectionTPclass} we introduce the ``triple-point''
classes,
or Gromov-Witten classes,
from which the quantum product is defined.
Gromov-Witten classes were first rigorously defined by Ruan
\cite{ruan} using symplectic deformations.
We have found it convenient to define these classes
(which measure rational curves on $\surf$ with $3$ marked points
satisfying certain geometric conditions)
in terms of the locus of such curves inside linear systems on
$\surf$,
rather than using the space of maps from ${\Bbb P}^1$ to $\surf$.
Other definitions (see
\cite{kontsevich-manin,ruan,siebert-tian,witten1,witten2})
use the space of maps instead; for rational surfaces this seems
unnecessary, at least for the definitions,
and slightly less transparent for our construction.
In Section \ref{sectiondefqntm} we give the definition
of the quantum cohomology ring.

As noted above, the quantum product is defined in terms of classes
of loci of rational curves on $\surf$; not all such loci
or cohomology classes appear in the definition,
and we call those that do ``relevant'' classes.
In Section \ref{strictrelevantsection} we enumerate the relevant
classes
for the general strict Del Pezzo surfaces,
that is, for the general blow-up of the plane at $6$ or fewer
points.
(These are the surfaces which are Fano surfaces in the strictest sense,
namely that $-K$ is very ample.)
In Section \ref{sectionformula} we give explicit formulas
for the quantum product in terms of the relevant classes on
$\surf$,
and using these formulas in Section \ref{sectionordinary} we show
that there is a natural map from the quantum cohomology ring
to the ordinary cohomology ring.
In Section \ref{sectionexamples} we explicitly compute
the quantum cohomology rings using the formulas
for the minimal rational strict Del Pezzo surfaces,
namely ${\Bbb P}^2$, ${\Bbb F}_0 = {\Bbb P}^1\times{\Bbb P}^1$,
and ${\Bbb F}_1$.

Quantum cohomology seems to be an admixture of homology and
cohomology;
as such, a functoriality property seems elusive
(e.g., should it be contravariant or covariant?).
We show in Section \ref{sectionfunctor}
that there is at least some degree of functoriality
for a blowup.  We use this to give a new proof of the associativity
of the quantum product for the general strict Del Pezzo surface
in Section \ref{sectionSDPassoc}; the functoriality property is
strong
enough so that it suffices to check associativity on the general
$6$-fold
blowup of the plane.
Moreover the results hold without the genericity assumptions on the
$6$-fold blowup mentioned above, although it is known in this case
that the general $n$-fold blowup of the plane is general in the
above sense
for $n \leq 9$.
We view this section as the primary technical contribution
of the paper, namely the verification of the associativity of the
quantum
product using algebro-geometric techniques rather than methods from
symplectic geometry.

Associativity in general is a tricky business;
the original heuristic arguments from physics
have been made precise only using symplectic geometry
(\cite{mcduff-salamon,ruan-tian}). We have given in Section
\ref{sectionassoc}
an alternate approach using algebraic geometry.
The sketch which we offer here is related to the argument
given in \cite{witten2}.
Finally in Section \ref{sectionenum} we draw some of the
enumerative
consequences of the associativity of the quantum product.
These lead in particular to the formulas appearing in
\cite{kontsevich-manin}, interpreted properly.
As a small application we show that associativity for the general
cubic
surface is essentially equivalent to the standard enumerative facts
concerning lines: there are $27$ of them and each meets $10$
others.

The authors would like to thank Igor Dolgachev,
who inspired the second author to think about these questions
with an excellent lecture on the subject.  We also profited
greatly from conversations with Antonella Grassi, Sheldon Katz,
and David Morrison.

\section{Classes on the triple product}
\label{sectionX3class}

Let $\surf$ be a general rational surface;
by this we mean firstly that
if $\surf$ is obtained by blowing up $n$ points $p_i$,
creating exceptional curves $E_i$, then for every $d$ and $m_i$,
the linear system $|dH - \sum_{i=1}^n m_i E_i|$ has the expected
dimension
(where $H$ is the pullback from the plane of the line class);
this expected dimension is ${1\over 2}[d^2+3d - \sum_i m_i(m_i+1)]$
unless this is negative.
Secondly we also assume that the locus of irreducible rational
curves in this
system has codimension equal to the arithmetic genus of the general
curve
in the system, which is the expected codimension.

The ordinary integral cohomology $H^*(\surf) = H^*(\surf,{\Bbb Z})$
is
\begin{eqnarray*}
H^0(\surf) &=& {\Bbb Z}[\surf], \\
H^1(\surf) &=& \{0\}, \\
H^2(\surf) &\cong& {\Bbb Z}^\rho, \\
H^3(\surf) &=& \{0\}, \\
H^4(\surf) &=& {\Bbb Z}[p], \text{ and } \\
H^i(\surf) &=& \{0\} \text{ for } i \geq 5,
\end{eqnarray*}
where here $[\surf]$ is the fundamental class of the surface
$\surf$ itself,
and $[p]$ is the class of a single point.
The Picard number $\rho$ is the rank of the $H^2$ group.

The triple product $\surf^3 = \surf \times \surf \times \surf$
is then an algebraic six-fold.
Its cohomology is (by the K\"unneth Theorem)
the triple tensor product of the cohomology of $\surf$,
and is therefore generated over ${\Bbb Z}$ by tensors of the form
$\alpha \otimes \beta \otimes \gamma$,
where $\alpha$, $\beta$, and $\gamma$ are either $[\surf]$, $[p]$,
or generators for $H^2(\surf,{\Bbb Z})$.
In particular all cohomology is even-dimensional.

Suppose that $[A]$ is a cohomology class in the triple product,
so that $[A] \in H^{2d}(\surf^3)$,
where $d$ is the complex codimension of the class.
Suppose that we choose $\alpha$ and $\beta$ to be homogeneous
elements
of $H^*(\surf)$, of degrees $2a$ and $2b$,
such that
\begin{equation}
\label{dimcondition}
4 \leq a + b + d \leq 6.
\end{equation}
In this case if we let $c = 6 - a - b - d$,
then for any class $\gamma \in H^{2c}(\surf)$,
the class $(\alpha\otimes\beta\otimes\gamma)$
will have the complementary dimension to the class $[A]$,
and therefore the intersection product
\[
[A] \cdot (\alpha \otimes \beta \otimes \gamma) \in {\Bbb Z}
\]
will be defined.
We therefore obtain a linear functional
\[
\Phi_{[A]}(\alpha,\beta): H^{2c}(\surf) \to {\Bbb Z}
\]
which by duality must be represented
by intersection with a cohomology class in $H^{4-2c}(\surf)$.
Call this cohomology class $\phi_{[A]}(\alpha,\beta)$;
in this case we have by definition
that for all $\gamma \in H^{2c}(\surf)$,
\[
[A] \cdot (\alpha \otimes \beta \otimes \gamma)
= \phi_{[A]}(\alpha,\beta) \cdot \gamma,
\]
and indeed by duality this characterizes the class
$\phi_{[A]}(\alpha,\beta)$.
Note that the intersection on the left side of this formula
is intersection in the cohomology of the triple product $\surf^3$,
while the intersection on the right side is intersection
in the cohomology of $\surf$.

The element $\phi_{[A]}(\alpha,\beta)$, by definition,
is linear in both $\alpha$ and $\beta$.
For notational sanity we declare $\phi_{[A]}(\alpha,\beta) = 0$
unless we have the dimension condition that $4 \leq a + b + d \leq
6$.

\begin{example}
\label{phiDelta}
Let $[\Delta]$ be the class of the diagonal
$\Delta \subset \surf \times \surf \times \surf$.
Since $\Delta$ has complex codimension $d = 4$,
$[\Delta] \in H^8(\surf^3)$.
Let $a$, $b$, and $c$ be non-negative integers such that $a+b \leq
2$
(this is (\ref{dimcondition}))
and $c = 2 - a - b$.
In this case if $\alpha$, $\beta$, and $\gamma$ are classes in
$H^*(\surf)$
of degrees $a$, $b$, and $c$ respectively,
the intersection product $\alpha \cdot \beta \cdot \gamma$
is defined in ${\Bbb Z}$;
this is just cup product to $H^4(\surf)$, then taking the degree.
In particular it is equal to $(\alpha \cup \beta) \cdot \gamma$.
Moreover $[\Delta] \cdot \alpha\otimes\beta\otimes\gamma$
is equal to this triple intersection.
Hence we see that $\phi_{[\Delta]}(\alpha,\beta) =
\alpha\cup\beta$.
\end{example}

\begin{example}
Let $[A]$ be the fundamental class of $\surf^3$;
this has codimension $d = 0$,
and lies in $H^0(\surf^3)$.
Indeed, it is equal to $[\surf]\otimes[\surf]\otimes[\surf]$.
Suppose that $a$, $b$, and $c=6-a-b$ are possible
complex codimensions in $\surf$,
satisfying (\ref{dimcondition}),
which is that $4 \leq a+b \leq 6$;
since $\surf$ is a surface, we must have $a = b = c = 2$.
If $\alpha$, $\beta$, and $\gamma$ are classes in $H^4(\surf)$,
by linearity for the computation we may take all three classes
equal to the
class $[p]$ of a point.
In this case we obviously have $[A] \cdot [p]\otimes[p]\otimes[p]
= 1$.
Therefore $\phi_{[A]}([p],[p]) = [\surf]$,
and unless both $\alpha$ and $\beta$ lie in $H^4(\surf)$,
$\phi_{[A]}(\alpha,\beta) = 0$.
\end{example}

\begin{example}
\label{exCCC}
Let $C$ be an divisor on $\surf$,
and let $[A] = [C] \otimes [C] \otimes [C] \in H^6(\surf^3)$.
Suppose that $a$, $b$, and $c=6-a-b$ are possible
complex codimensions in $\surf$,
satisfying (\ref{dimcondition}),
which is that $1 \leq a+b \leq 3$,
and $\alpha$, $\beta$, and $\gamma$ are classes in $H^*(\surf)$
of degrees $a$, $b$, and $c$ respectively.
Then $[A] \cdot \alpha\otimes\beta\otimes\gamma =
([C]\cdot \alpha)([C]\cdot \beta)([C]\cdot \gamma)$,
which is zero unless $a=b=c=1$,
and $\alpha$, $\beta$, $\gamma$ are divisor classes.
Therefore in this case
$\phi_{[A]}(\alpha,\beta) = ([C]\cdot \alpha)([C]\cdot \beta)[C]$
for divisor classes $\alpha$ and $\beta$,
and is zero otherwise.
\end{example}

\begin{example}
Let $C$ be an divisor on $\surf$,
and let $[A] = [C] \otimes [C] \otimes [\surf] \in H^4(\surf^3)$.
Suppose that $\alpha$ and $\beta$ are divisor classes in
$H^2(\surf)$.
Then $[A] \cdot \alpha\otimes\beta\otimes[p] =
([C]\cdot \alpha)([C]\cdot \beta)$.
Therefore in this case
$\phi_{[A]}(\alpha,\beta) = ([C]\cdot \alpha)([C]\cdot
\beta)[\surf]$
for divisor classes $\alpha$ and $\beta$,
and is zero otherwise.
\end{example}

\begin{example}
Let $C$ be an divisor on $\surf$,
and instead let $[A] = [C] \otimes [\surf] \otimes [C] \in
H^4(\surf^3)$.
Suppose that $\alpha$ and $\gamma$ are divisor classes in
$H^2(\surf)$.
Then $[A] \cdot \alpha\otimes[p]\otimes\gamma =
([C]\cdot \alpha)([C]\cdot \gamma)$.
Therefore in this case
$\phi_{[A]}(\alpha,[p]) = ([C]\cdot \alpha)[C]$
for a divisor class $\alpha$
(and is zero otherwise).
\end{example}

\begin{example}
Let $C$ be an divisor on $\surf$,
and let $[A] = [C] \otimes [\surf] \otimes [\surf] \in
H^2(\surf^3)$.
Suppose that $\alpha$ is a divisor classes in $H^2(\surf)$.
Then $[A] \cdot \alpha\otimes[p]\otimes[p] =
([C]\cdot \alpha)$.
Therefore in this case
$\phi_{[A]}(\alpha,[p]) = ([C]\cdot \alpha)[\surf]$
for a divisor class $\alpha$
(and is zero otherwise).
\end{example}

\begin{example}
Let $C$ be an divisor on $\surf$,
and instead let $[A] = [\surf] \otimes [\surf] \otimes [C] \in
H^2(\surf^3)$.
Then for any divisor class $\gamma$ in $H^2(\surf)$,
we have $[A] \cdot [p]\otimes[p]\otimes\gamma =
([C]\cdot \gamma)$.
Therefore in this case
$\phi_{[A]}([p],[p]) = [C]$.
\end{example}

\section{Three-point classes on the triple product}
\label{sectionTPclass}

Fix a divisor class $[C] \in H^2(\surf)$,
such that $|C|$ has no fixed components and is non-empty.
Define the locus ${\cal R}_{[C]} \subset |C|$
representing irreducible rational curves with only nodes as
singularities.
Inside the product ${\cal R}_{[C]} \times \surf^3$
consider the incidence correspondence
\[
{\cal S}_{[C]} = \{(C,p_1,p_2,p_3)\;|\; \text{ the points }p_i
\text{ are distinct smooth points on the curve }C\}.
\]
Let $\overline{{\cal S}}_{[C]}$ be the closure of this subvariety
inside $|C|\times \surf^3$.
The second projection gives a regular map
$\pi:\overline{{\cal S}}_{[C]} \to \surf^3$.
We define the {\em three-point class} $[A_{[C]}]$
to be $\pi_*[\overline{{\cal S}}_{[C]}]$,
the image of the fundamental class.
(This is a priori in homology,
but we consider it in cohomology
via duality.)

Note that $\dim {\cal S}_{[C]} = 3 + \dim {\cal R}_{[C]}$.

A degenerate version of this locus is obtained
when we allow the cohomology class $[C]\in H^2(\surf)$ to be
trivial;
we declare in this case that the three-point class $[A_{[0]}]$ is
the class
of the diagonal $\Delta$.

It is an elementary matter to compute the dimensions
of these three-point classes,
in terms of the numerical characters of the class $[C]$.
Since $\surf$ is a general rational surface,
the general member of the linear system $|C|$
is smooth, and the system has the  expected dimension,
which is $(C\cdot C) +1-p_a(C)$;
here $p_a(C)$ is the arithmetic genus, and equals
$p_a(C) = 1 + (C \cdot C + K_{\surf})/2$ by Riemann-Roch.

Imposing a node on a member of $|C|$ is one condition;
hence the locus of nodal rational curves ${\cal R}_{[C]}$
has dimension $\dim |C| - p_a(C) = (C \cdot C) + 1 - 2p_a(C)$
by our general assumption on $\surf$.
Therefore the dimension of the incidence locus ${\cal S}_{[C]}$
is $(C \cdot C) + 4 - 2p_a(C)$,
which we may re-write as
\[
\dim {\cal S}_{[C]} = 2 - (C \cdot K_{\surf}).
\]

To be more explicit,
suppose that $[C] = dH - \sum_{i=1}^n m_i E_i$,
where $H$ is the pullback of the line class from ${\Bbb P}^2$
and $E_i$ is the class of the exceptional curve over the blown up
point $p_i$.
Then $K_{\surf} = -3H + \sum_i E_i$,
so that $(C \cdot K_{\surf}) = -3d + \sum_i m_i$.
Hence we have
\[
\dim {\cal S}_{[C]} = 3d +2 - \sum_i m_i.
\]

The only classes of curves on $\surf$ which are not of this form
are the classes of the exceptional curves $E_i$ themselves.
Here $|E_i|$ is a single point (the only member is $E_i$ itself)
and $E_i$ is a smooth rational curve;
so $\dim {\cal R}_{[E_i]} = 0$ and $\dim {\cal S}_{[E_i]} = 3$.
(Actually the formula holds in this case also,
with $d=0$, $m_i = -1$, and $m_j = 0$ for $j \neq i$.)

\begin{definition}
A class $[C] \in H^2(\surf)$ is {\em relevant} (for quantum
cohomology)
if either $[C] = 0$ or ${\cal R}_{[C]} \neq \emptyset$
and $\dim {\cal S}_{[C]} \leq 6$.
\end{definition}

If $\dim {\cal S}_{[C]} > 6$, then the image of the fundamental
class
of its closure in the six-dimensional variety $\surf^3$
is trivial.  Therefore all non-relevant classes induce a trivial
three-point class $[A_{[C]}]$.

Given the class $[A_{[C]}]$, they induce as noted in the previous
section
classes $\phi_{[A_{[C]}]}(\alpha,\beta)$ in the cohomology of
$\surf$.
We will abbreviate the notation for these classes and write simply
$\phi_{[C]}(\alpha,\beta)$.

We note that there is an obvious $S_3$-symmetry to the three-point
classes,
in the sense that
\[
[A_{[C]}] \cdot \alpha_1 \otimes \alpha_2 \otimes \alpha_3
=
[A_{[C]}] \cdot
\alpha_{\sigma(1)} \otimes \alpha_{\sigma(1)} \otimes
\alpha_{\sigma(1)}
\]
for any permutation $\sigma \in S_3$.
This is simply because the locus ${\cal S}_{[C]}$ is
$S_3$-invariant.
As a consequence of this, we see that the $\phi$-classes are
symmetric:
\[
\phi_{[C]}(\alpha,\beta) = \phi_{[C]}(\beta,\alpha)
\]
and of course they are bilinear in $\alpha$ and $\beta$.

\begin{example}
\label{phi0}
If we start with the trivial class $[C] = 0$,
then the three-point class $[A_{[0]}]$ is the class of the
diagonal.
Hence as we have noted above,
for any classes $\alpha$ and $\beta$ in $H^*(\surf)$,
\[
\phi_0(\alpha,\beta) = \alpha \cup \beta.
\]
\end{example}

\begin{example}
Suppose that $E$ is a $(-1)$-curve on $\surf$,
that is, a smooth rational curve with $(E \cdot E) = -1$.
The only member of the linear system $|E|$ is the curve $E$ itself;
the locus ${\cal R}_{[E]}$ of nodal rational curves in $|E|$
is the whole system $|E| = \{E\}$.
The incidence locus ${\cal S}_{[E]} = \{E\} \times E \times E
\times E$
with the large diagonal removed;
its closure $\overline{{\cal S}}_{[E]} = \{E\} \times E \times E
\times E$.
Hence the image under the projection to $\surf^3$
is $E\times E\times E$, and the class $[A_{[E]}] =
[E]\otimes[E]\otimes[E]$.
Hence by the computation in Example \ref{exCCC},
we have
\[
\phi_{[E]}(\alpha,\beta) = ([E]\cdot \alpha)([E]\cdot \beta)[E]
\]
for divisor classes $\alpha$ and $\beta$,
and is zero unless $\alpha$ and $\beta$ are both in $H^2(M)$.
\end{example}

\begin{example}
More generally suppose that $E$ is an irreducible curve on $\surf$
whose class is relevant, and $(E\cdot E) - 2p_a(E) = -1$,
so that $\dim {\cal R}_{[E]} = 0$ and is therefore a finite set;
say it has $d$ members $E_1,\dots,E_d$.
The incidence locus ${\cal S}_{[E]} =
\bigcup_{i=1}^d \{E_i\} \times E_i \times E_i \times E_i$
with the large diagonal removed;
its closure $\overline{{\cal S}}_{[E]} =
\bigcup_{i=1}^d \{E_i\} \times E_i \times E_i \times E_i$.
Hence the image under the projection to $\surf^3$
is $\bigcup_{i=1}^d E_i\times E_i\times E_i$,
and the class $[A_{[E]}] = d [E]\otimes[E]\otimes[E]$.
Hence by the computation in Example \ref{exCCC},
we have
\[
\phi_{[E]}(\alpha,\beta) = d ([E]\cdot \alpha)([E]\cdot \beta)[E]
\]
for divisor classes $\alpha$ and $\beta$,
and is zero unless $\alpha$ and $\beta$ are both in $H^2(M)$.
This generalizes the previous example, where $d=1$.
\end{example}

\begin{example}
Suppose that $F$ is a fiber in a ruling on $\surf$,
that is, a smooth rational curve with $(F \cdot F) = 0$.
The linear system $|F|$ is a pencil;
the locus ${\cal R}_{[F]}$ of irreducible nodal rational curves in
$|F|$
is an open dense subset of the whole system $|F|$
(it is the set of smooth members of $|F|$).
The incidence locus ${\cal S}_{[F]}$ has complex dimension $4$;
an element is obtained by choosing a member of $|F|$,
then three points on this member.
The complementary classes in $H^*(\surf^3)$ have complex
codimension $4$,
that is, they are the classes in $H^8(\surf^3)$.
This group is generated by the classes
$[p]\otimes[p]\otimes[\surf]$,
$\alpha\otimes\beta\otimes[p]$ for divisor classes $\alpha$ and
$\beta$,
and the associated classes obtained by symmetry.
The intersection product
$[A_{[F]}] \cdot [p]\otimes[p]\otimes[\surf] = 0$,
since there is a no curve in the system through two general points
$p$.
The intersection product
$[A_{[F]}] \cdot \alpha\otimes\beta\otimes[p] =
(F\cdot \alpha)(F \cdot \beta)$;
forcing the curve to pass through the general point $p$ gives a
unique
curve in the system, and the choice of the other two points,
which must lie in the divisor $\alpha$ and $\beta$ respectively,
gives the result above.
Therefore we have
\[
\phi_{[F]}(\alpha,\beta) = ([F]\cdot \alpha)([F]\cdot \beta)[\surf]
\]
for divisor classes $\alpha$ and $\beta$.  Moreover by the symmetry
we also have
$[A_{[F]}] \cdot \alpha\otimes[p]\otimes\gamma =
(F\cdot \alpha)(F \cdot \gamma)$, so that (after taking symmetry
into account)
\[
\phi_{[F]}([p],\alpha) = \phi_{[F]}(\alpha,[p]) = ([F]\cdot
\alpha)[F]
\]
for a divisor class $\alpha$.  All other $\phi$-classes are zero.
\end{example}

\begin{example}
Suppose that $F$ gives a relevant class on $\surf$
with $(F\cdot F) - 2p_a(F) = 0$;
then the locus ${\cal R}_{[F]}$ of nodal rational curves in $|F|$
forms a curve.
Denote by $d$ the degree of this curve in the projective space
$|F|$.
The incidence locus ${\cal S}_{[F]}$ has complex dimension $4$;
an element is obtained by choosing a member of ${\cal R}_{[F]}$,
then three points on this member.
The complementary classes are the classes in $H^8(\surf^3)$;
This group as above is generated by the classes
$[p]\otimes[p]\otimes[\surf]$,
$\alpha\otimes\beta\otimes[p]$ for divisor classes $\alpha$ and
$\beta$,
and the associated classes obtained by symmetry.
The intersection product
$[A_{[F]}] \cdot [p]\otimes[p]\otimes[\surf] = 0$,
since there is a no curve in the system through two general points
$p$.
The intersection product
$[A_{[F]}] \cdot \alpha\otimes\beta\otimes[p] =
d (F\cdot \alpha)(F \cdot \beta)$;
forcing the curve to pass through the general point $p$
gives $d$ curves in the system,
and the choice of the other two points,
which must lie in the divisor $\alpha$ and $\beta$ respectively,
contributes $(F\cdot \alpha)$ and $(F\cdot \beta)$ respectively
to the number of choices.
Therefore we have
\[
\phi_{[F]}(\alpha,\beta) = d ([F]\cdot \alpha)([F]\cdot
\beta)[\surf]
\]
for divisor classes $\alpha$ and $\beta$.  Moreover by the symmetry
we also have
$[A_{[F]}] \cdot \alpha\otimes[p]\otimes\gamma =
d (F\cdot \alpha)(F \cdot \gamma)$,
so that (after taking symmetry into account)
\[
\phi_{[F]}([p],\alpha) = \phi_{[F]}(\alpha,[p]) = d([F]\cdot
\alpha)[F]
\]
for a divisor class $\alpha$.  All other $\phi$-classes are zero.
This generalizes the previous example, where $d=1$.
\end{example}

\begin{example}
Suppose that $L$ is a smooth rational curve on $\surf$
with $(L \cdot L) = 1$.
The linear system $|L|$ is a net;
the locus ${\cal R}_{[L]}$ of nodal (i.e. smooth) rational curves
in $|L|$
is an open dense subset of $|L|$.
The incidence locus ${\cal S}_{[L]}$ has complex dimension $5$;
an element is obtained by choosing a smooth member of $|L|$,
then three points on this member.
The complementary classes in $H^*(\surf^3)$ have complex
codimension $5$,
that is, they are the classes in $H^{10}(\surf^3)$.
This group is generated by the classes
$[p]\otimes[p]\otimes\alpha$,
for a divisor classes $\alpha$,
and the associated classes obtained by symmetry.
The intersection product
$[A_{[L]}] \cdot [p]\otimes[p]\otimes\alpha = (L \cdot \alpha)$,
since through two general points there is a unique member $L_0$ of
$|L|$,
whose third point can be any of the points
where $L_0$ meets the divisor $\alpha$.
Therefore we have
\[
\phi_{[L]}(\alpha,[p]) = (L\cdot \alpha)[\surf]
\]
for a divisor classes $\alpha$, and
\[
\phi_{[L]}([p],[p]) = [L].
\]
All other $\phi$-classes are zero.
\end{example}

\begin{example}
Again we may generalize the above in case
$L$ induces any relevant class with $(L \cdot L) -2 p_a(L) = 1$.
The locus ${\cal R}_{[L]}$ of nodal (i.e. smooth) rational curves
in $|L|$
is a surface inside the projective space $|L|$;
let $d$ be the degree of this surface.
The incidence locus ${\cal S}_{[L]}$ has complex dimension $5$.
The complementary classes in $H^*(\surf^3)$
are the classes in $H^{10}(\surf^3)$,
which is generated by the classes
$[p]\otimes[p]\otimes\alpha$,
for a divisor classes $\alpha$,
and the associated classes obtained by symmetry.
The intersection product
$[A_{[L]}] \cdot [p]\otimes[p]\otimes\alpha = d (L \cdot \alpha)$,
since through two general points there are $d$ members of ${\cal
R}_{[L]}$,
whose third point can be any of the points
where the member meets the divisor $\alpha$.
Therefore we have
\[
\phi_{[L]}(\alpha,[p]) = d (L\cdot \alpha)[\surf]
\]
for a divisor classes $\alpha$, and
\[
\phi_{[L]}([p],[p]) = d [L].
\]
All other $\phi$-classes are zero.
This generalizes the previous example, where $d=1$.
\end{example}

\begin{example}
Suppose that $C$ is a smooth rational curve on $\surf$
with $(C \cdot C) = 2$.
The linear system $|C|$ is a web (that is, it is $3$-dimensional);
the locus ${\cal R}_{[C]}$ of nodal (i.e. smooth) rational curves
in $|C|$
is again an open dense subset of $|C|$.
The incidence locus ${\cal S}_{[C]}$ has complex dimension $6$;
an element is obtained by choosing a smooth member of $|C|$,
then three points on this member.
The complementary classes in $H^*(\surf^3)$ have complex
codimension $6$,
that is, they are the classes in $H^{12}(\surf^3)$.
This group has rank one, and is generated by the class
$[p]\otimes[p]\otimes[p]$.
The intersection product
$[A_{[C]}] \cdot [p]\otimes[p]\otimes[p] = 1$,
since through three general points there is a unique member of
$|C|$.
Indeed, we have that $[A_{[C]}] =
[\surf]\otimes[\surf]\otimes[\surf]$.
Therefore we have
\[
\phi_{[C]}([p],[p]) = [\surf].
\]
All other $\phi$-classes are zero.
\end{example}

\begin{example}
\label{phiC2=2}
Again if $C$ is a relevant class with $(C\cdot C) - 2p_a(C) = 2$,
then the locus ${\cal R}_{[C]}$ of nodal rational curves in $|C|$
is a threefold of $|C|$;
let $d$ be the degree of this threefold.
The incidence locus ${\cal S}_{[C]}$ has complex dimension $6$;
and the complementary classes in $H^*(\surf^3)$ have complex
codimension $6$,
that is, they are the classes in $H^{12}(\surf^3)$,
which is generated by the class
$[p]\otimes[p]\otimes[p]$.
The intersection product
$[A_{[C]}] \cdot [p]\otimes[p]\otimes[p] = d$,
since through three general points
there is are $d$ members of ${\cal R}_{[C]}$.
Indeed, we have that $[A_{[C]}] = d
[\surf]\otimes[\surf]\otimes[\surf]$.
Therefore we have
\[
\phi_{[C]}([p],[p]) = d [\surf].
\]
All other $\phi$-classes are zero.
This generalizes the previous example, where $d=1$.
\end{example}

We offer the following example which shows that the above
phenomenon
occurs, namely that there are relevant classes
which come from singular rational curves.

\begin{example}
Let $\surf$ be the blow-up of the plane ${\Bbb P}^2$
at $5$ general points $p_1,\dots p_5$.
Let $H$ denote the line class on $\surf$
and $E_i$ denote the exceptional curve lying over $p_i$.
Consider the anti-canonical class $C = -K_{\surf} = 3H -
\sum_{i=1}^5 E_i$.
This is the linear system of cubics passing through the $5$ points
$p_i$.
Note that $(C\cdot C) = 4$ and $p_a(C) = 1$
so that $(C\cdot C) - 2p_a(C) = 2$.
We have $\dim {\cal S}_{[C]} = 2 - (C \cdot K_{\surf}) = 6$,
so $[C]$ is a relevant class.
The map
\[
\pi_2:\overline{{\cal S}_{[C]}} \to \surf^3
\]
has as its general fiber over a triple $(q_1,q_2,q_3)$
those nodal rational curves in the linear system $|C|$
through the three points $q_1$, $q_2$, and $q_3$.
This is exactly the set of nodal rational cubics in the plane
passing through the eight points $p_1,\dots p_5$, $q_1,\dots q_3$.
The system of cubics through these general eight points
forms a pencil of genus one curves,
which has exactly $12$ singular members.
Hence the map $\pi_2:\overline{{\cal S}_{[C]}} \to \surf^3$
is generically finite of degree $12$,
and so pushing down the fundamental class we see that
\[
[A_{[C]}] = 12[\surf]\otimes[\surf]\otimes[\surf].
\]
Hence
\[
[A_{[C]}] \cdot [p]\otimes[p]\otimes[p] = 12 \text{ and }
\phi_{[C]}([p],[p]) = 12[\surf].
\]
All other $\phi$-classes are zero.
\end{example}

\section{Definition of quantum cohomology}
\label{sectiondefqntm}

We denote by $\Eff(\surf)$
the cone of effective divisor classes in $H^2(\surf)$.
It forms a semigroup under addition.

Let $Q = {\Bbb Z}[[\Eff(\surf)]]$
be the completion of the integral group ring over $\Eff(\surf)$;
$Q$ is a ${\Bbb Z}$-algebra.
It is customary to formally introduce a multi-variable $q$
and to write the module generators of $Q$ as elements $q^{[D]}$,
where $[D]$ is an effective cohomology class in $\Eff(\surf)$.
With this notation, every element of $Q$ can be written
as a formal series
\[
\sum_{[D]\in\Eff(\surf)} n_{[D]} q^{[D]}
\]
with integral coefficients $n_{[D]}$,
and divisor class exponents $[D] \in \Eff(\surf)$.
In this way multiplication in the ring $Q$
is induced by the relations that
\[
q^{[D_1]}q^{[D_2]} = q^{[D_1+D_2]}
\]
for divisors $D_1$ and $D_2$ on $\surf$,
and the ordinary distributive and associative laws.

Define the {\em quantum cohomology ring} of $\surf$ to be
\[
H^*_Q(\surf) = H^*(\surf) \otimes_{\Bbb Z} Q
\]
as a free abelian group.
Moreover it is also a $Q$-module, with the obvious structure.

The multiplication $\qprod$ on $H^*_Q(\surf)$,
called the {\em quantum product},
is determined by knowing the products of (homogeneous) elements
from
$H^*(\surf)$,
since the rest comes from linearity and the $Q$-module structure.
For two homogeneous elements $\alpha$ and $\beta$ in $H^*(\surf)$
define
\[
\alpha \qprod \beta = \sum_{[C]} \phi_{[C]}(\alpha,\beta) q^{[C]},
\]
the sum begin taken over
the relevant cohomology classes in $H^2(\surf)$.

We remark that if there are only
finitely many relevant classes in $H^2(\surf)$,
then the quantum cohomology ring may be formulated as a polynomial
ring
instead of a power series ring;
in other words, one may take $Q = {\Bbb Z}[\Eff(\surf)]$
to be simply the integral semigroup ring instead of its completion.
This is the case for a Del Pezzo surface $\surf$.

We have an immediate identification for the identity of the quantum
product:

\begin{lemma}
The fundamental class $[\surf]$ is the identity for the quantum
product.
In other words, for every class $\alpha \in H^*(\surf)$,
\[
\phi_0([\surf],\alpha) = \alpha
\]
and if $[C] \neq 0$, then
\[
\phi_{[C]}([\surf],\alpha) = 0.
\]
\end{lemma}

\begin{pf}
The $[C]= 0$ statement follows from the computation in Example
\ref{phiDelta};
we have
\[
\phi_0([\surf],\alpha) = \phi_{[\Delta]}([\surf],\alpha)
= [\surf]\cup\alpha = \alpha
\]
since $[\surf]$ is the identity for the ordinary cup product.
If $[C] \neq 0$, then
$\dim {\cal S}_{[C]} = 3 + \dim {\cal R}_{[C]} \geq 3$ if $[C]$ is
relevant.
Therefore the three-point class $[A_{[C]}]$
lies in $H^{2k}(\surf^3)$ for $k \leq 3$.
Any complementary class of the form
$[\surf]\otimes\beta\otimes\gamma$
must have $\deg(\beta)+\deg(\gamma) = 12-2k \leq 6$.

Suppose first that $\dim {\cal S}_{[C]} = 3$,
so that ${\cal R}_{[C]}$ is a finite set and
the three-point class $[A_{[C]}]$
lies in $H^{6}(\surf^3)$.
The only complementary classes involving the fundamental class
$[\surf]$
have the form $[\surf]\otimes\beta\otimes[p]$ for some divisor
class $\beta$.
But $[A_{[C]}] \cdot [\surf]\otimes\beta\otimes[p]$ counts the
number of curves
in ${\cal S}_{[C]}$ whose second point lies in the divisor $\beta$
and whose third point equals $p$;
since $p$ is a general point and ${\cal R}_{[C]}$ is a finite set,
there are no curves in ${\cal R}_{[C]}$ through $p$
and this intersection number is zero.
Hence $\phi_{[C]}([\surf],\alpha) = 0$ for any $\alpha$.

Suppose next that $\dim {\cal S}_{[C]} = 4$,
so that the locus ${\cal R}_{[C]}$ is one-dimensional and
the three-point class $[A_{[C]}]$
lies in $H^{4}(\surf^3)$.
The only complementary classes involving the fundamental class
$[\surf]$
have the form $[\surf]\otimes[p]\otimes[p]$.
The intersection product
$[A_{[C]}] \cdot [\surf]\otimes[p]\otimes[p]$ counts the number of
curves
in ${\cal S}_{[C]}$ whose second and third point are specified
general points;
since ${\cal R}_{[C]}$ is one-dimensional,
there are no curves in ${\cal R}_{[C]}$ through two specified
general points,
and this intersection number is zero.
Hence again $\phi_{[C]}([\surf],\alpha) = 0$ for any $\alpha$.

Finally if $\dim {\cal S}_{[C]} \geq 5$,
there are no complementary classes in the cohomology of $\surf^3$
of the form $[\surf]\otimes\beta\otimes\gamma$ at all.
\end{pf}

\section{Relevant classes on strict Del Pezzo surfaces}
\label{strictrelevantsection}

Let $\surf$ be a general strict Del Pezzo surface,
that is, $\surf = {\Bbb F}_0 \cong {\Bbb P}^1\times {\Bbb P}^1$
or a general blowup of the plane such that $-K_{\surf}$ is very
ample.
This amounts to having $X \cong {\Bbb F}_0$
or $X$ being a blowup of the plane at $n \leq 6$ general points.
It is an elementary matter to compute the relevant classes
on such a strict Del Pezzo surface $X$.
The results are shown in the tables below.

In the first few columns are the numerical characters of the class:
the bidegree in the case of ${\Bbb F}_0$
and the integers $d$ and $m_i$ for a class of the form $dH - \sum_i
m_i E_i$
on a blowup of the plane; here $d$ is the degree and $m_i$
is the multiplicity of the curves at the associated blown up point.
The final three columns contain the quantity $C^2-2p_a(C)$
(on which relevance is based), the arithmetic genus $p_a(C)$,
and the number of such classes up to permutations of the $E_i$'s.

\begin{center}
\begin{tabular}{c|c|c|c}
\multicolumn{4}{c}{Relevant nonzero classes on ${\Bbb P}^2$} \\
\hline
degree & $C^2-2p_a(C)$ & $p_a(C)$ & $\#$ \\
1 & 1 & 0 & 1 \\
\end{tabular}
\end{center}

\begin{center}
\begin{tabular}{c|c|c|c}
\multicolumn{4}{c}{Relevant classes on ${\Bbb F}_0$} \\ \hline
bidegree & $C^2-2p_a(C)$ & $p_a(C)$ &  $\#$ \\
(0,1) & 0 & 0 & 1 \\
(1,0) & 0 & 0 & 1 \\
(1,1) & 2 & 0 & 1 \\
\end{tabular}
\end{center}

Denote by $X_n$ a general blowup of ${\Bbb P}^2$ at $n$ points.

\begin{center}
\begin{tabular}{c|c|c|c|c}
\multicolumn{5}{c}{Relevant classes on $X_1$} \\ \hline
degree & $m_1$ & $C^2-2p_a(C)$ & $p_a(C)$ & $\#$ \\
0 & -1 & -1 & 0 & 1 \\
1 &  1 &  0 & 0 & 1 \\
1 &  0 &  1 & 0 & 1 \\
\end{tabular}
\end{center}

\begin{center}
\begin{tabular}{c|cc|c|c|c}
\multicolumn{6}{c}{Relevant nonzero classes on $X_2$} \\ \hline
degree & $m_1$ & $m_2$ & $C^2-2p_a(C)$ & $p_a(C)$ & $\#$ \\
0 & -1 & 0 & -1 & 0 & 2 \\
1 &  1 & 1 & -1 & 0 & 1 \\
1 &  1 & 0 &  0 & 0 & 2 \\
1 &  0 & 0 &  1 & 0 & 1 \\
2 &  1 & 1 &  2 & 0 & 1 \\
\end{tabular}
\end{center}

\begin{center}
\begin{tabular}{c|ccc|c|c|c}
\multicolumn{7}{c}{Relevant classes on $X_3$} \\ \hline
degree & $m_1$ & $m_2$ & $m_3$ & $C^2-2p_a(C)$ & $p_a(C)$ & $\#$ \\
0 & -1 & 0 & 0 & -1 & 0 & 3 \\
1 &  1 & 1 & 0 & -1 & 0 & 3 \\
1 &  1 & 0 & 0 &  0 & 0 & 3 \\
1 &  0 & 0 & 0 &  1 & 0 & 1 \\
2 &  1 & 1 & 1 &  1 & 0 & 1 \\
2 &  1 & 1 & 0 &  2 & 0 & 3 \\
\end{tabular}
\end{center}

\begin{center}
\begin{tabular}{c|cccc|c|c|c}
\multicolumn{8}{c}{Relevant nonzero classes on $X_4$} \\ \hline
degree & $m_1$ & $m_2$ & $m_3$ & $m_4$
& $C^2-2p_a(C)$ & $p_a(C)$ & $\#$ \\
0 & -1 & 0 & 0 & 0 & -1 & 0 & 4 \\
1 &  1 & 1 & 0 & 0 & -1 & 0 & 6 \\
1 &  1 & 0 & 0 & 0 &  0 & 0 & 4 \\
2 &  1 & 1 & 1 & 1 &  0 & 0 & 1 \\
1 &  0 & 0 & 0 & 0 &  1 & 0 & 1 \\
2 &  1 & 1 & 1 & 0 &  1 & 0 & 4 \\
2 &  1 & 1 & 0 & 0 &  2 & 0 & 6 \\
3 &  2 & 1 & 1 & 1 &  2 & 0 & 4 \\
\end{tabular}
\end{center}

\begin{center}
\begin{tabular}{c|ccccc|c|c|c}
\multicolumn{9}{c}{Relevant classes on $X_5$} \\ \hline
degree & $m_1$ & $m_2$ & $m_3$ & $m_4$ & $m_5$
& $C^2-2p_a(C)$ & $p_a(C)$ & $\#$ \\
0 & -1 & 0 & 0 & 0 & 0 & -1 & 0 & 5 \\
1 &  1 & 1 & 0 & 0 & 0 & -1 & 0 & 10 \\
2 &  1 & 1 & 1 & 1 & 1 & -1 & 0 & 1 \\
1 &  1 & 0 & 0 & 0 & 0 &  0 & 0 & 5 \\
2 &  1 & 1 & 1 & 1 & 0 &  0 & 0 & 5 \\
1 &  0 & 0 & 0 & 0 & 0 &  1 & 0 & 1 \\
2 &  1 & 1 & 1 & 0 & 0 &  1 & 0 & 10 \\
3 &  2 & 1 & 1 & 1 & 1 &  1 & 0 & 5 \\
2 &  1 & 1 & 0 & 0 & 0 &  2 & 0 & 10 \\
3 &  1 & 1 & 1 & 1 & 1 &  2 & 1 & 1 \\
3 &  2 & 1 & 1 & 1 & 0 &  2 & 0 & 20 \\
4 &  2 & 2 & 2 & 1 & 1 &  2 & 0 & 10 \\
\end{tabular}
\end{center}

\begin{center}
\begin{tabular}{c|cccccc|c|c|c}
\multicolumn{10}{c}{Relevant classes on $X_6$} \\ \hline
degree & $m_1$ & $m_2$ & $m_3$ & $m_4$ & $m_5$ & $m_6$
& $C^2-2p_a(C)$ & $p_a(C)$ & $\#$ \\
0 & -1 & 0 & 0 & 0 & 0 & 0 & -1 & 0 & 6 \\
1 &  1 & 1 & 0 & 0 & 0 & 0 & -1 & 0 & 15 \\
2 &  1 & 1 & 1 & 1 & 1 & 0 & -1 & 0 & 6 \\
1 &  1 & 0 & 0 & 0 & 0 & 0 &  0 & 0 & 6 \\
2 &  1 & 1 & 1 & 1 & 0 & 0 &  0 & 0 & 15 \\
3 &  2 & 1 & 1 & 1 & 1 & 1 &  0 & 0 & 6 \\
1 &  0 & 0 & 0 & 0 & 0 & 0 &  1 & 0 & 1 \\
2 &  1 & 1 & 1 & 0 & 0 & 0 &  1 & 0 & 20 \\
3 &  1 & 1 & 1 & 1 & 1 & 1 &  1 & 1 & 1 \\
3 &  2 & 1 & 1 & 1 & 1 & 0 &  1 & 0 & 30 \\
4 &  2 & 2 & 2 & 1 & 1 & 1 &  1 & 0 & 20 \\
5 &  2 & 2 & 2 & 2 & 2 & 2 &  1 & 0 & 1 \\
2 &  1 & 1 & 0 & 0 & 0 & 0 &  2 & 0 & 15 \\
3 &  1 & 1 & 1 & 1 & 1 & 0 &  2 & 1 & 6 \\
3 &  2 & 1 & 1 & 1 & 0 & 0 &  2 & 0 & 60 \\
4 &  2 & 2 & 1 & 1 & 1 & 1 &  2 & 1 & 15 \\
4 &  2 & 2 & 2 & 1 & 1 & 0 &  2 & 0 & 60 \\
4 &  3 & 1 & 1 & 1 & 1 & 1 &  2 & 0 & 6 \\
5 &  2 & 2 & 2 & 2 & 2 & 1 &  2 & 1 & 6 \\
5 &  3 & 2 & 2 & 2 & 1 & 1 &  2 & 0 & 60 \\
6 &  3 & 3 & 2 & 2 & 2 & 2 &  2 & 0 & 15 \\
\end{tabular}
\end{center}

The table of relevant classes on $X_6$ has some features
which will be useful below. Let us collect them in the following
lemma.
Note that the anti-canonical class $-K$ on the surface
has $d = 3$ and $m_i = 1$ for each $i = 1,\ldots,6$.

\begin{lemma}
\label{X6relevantlemma}
Let $X_6$ denote a general $6$-fold blowup of the plane
(that is, a general cubic surface in ${\Bbb P}^3$).
\begin{enumerate}
\item All relevant classes $[C]$ on $X_6$
have arithmetic genus $p_a(C) \leq 1$.
\item The anti-canonical class $-K$ is the unique relevant class
$[C]$ on $X_6$
with $C^2 - 2p_a(C) = 1$ and $p_a(C) = 1$.
\item There are exactly $27$ relevant classes $[E]$ on $X_6$
with $E^2 - 2p_a(E) = -1$; all have $p_a(E) = 0$.
These are the classes of the $27$ lines on the cubic surface.
\item There are exactly $27$ relevant classes $[F]$ on $X_6$
with $F^2 - 2p_a(F) = 0$; all have $p_a(F) = 0$.
Each such class $F$ is obtained from a unique relevant class $E$
with $E^2 - 2p_a(E) = -1$ by subtracting $E$ from the anticanonical
class:
$F = -K - E$.
Any two such classes $F$, $G$ satisfy $0 \leq (F \cdot G) \leq 2$;
$(F\cdot G) = 0$ if and only if $F = G$.
If $(F \cdot G) = 1$ then $C = F+G$ is a relevant class
with $C^2 - 2p_a(C) = 2$ and $p_a(C) = 0$.
If $(F \cdot G) = 2$ then $C = F+G$ is a relevant class
with $C^2 - 2p_a(C) = 2$ and $p_a(C) = 1$.
\item There are exactly $72$ relevant classes $[L]$ on $X_6$
with $L^2 - 2p_a(L) = 1$ and $p_a(L) = 0$.
For each such class $[L]$ and each relevant class $E$ with $E^2 -
2p_a(E) = -1$
we have $0 \leq (L \cdot E) \leq 2$.
If $(L\cdot E) = 1$ then $C = L+E$ is a relevant class
with $C^2 - 2p_a(C) = 2$ and $p_a(C) = 0$.
If $(L \cdot E) = 2$ then $C = L+E$ is a relevant class
with $C^2 - 2p_a(C) = 2$ and $p_a(C) = 1$.
\item There are exactly $216$ relevant classes $[C]$ on $X_6$
with $C^2 - 2p_a(C) = 2$ and $p_a(C) = 0$.
Each such class $C$ can be written uniquely (up to order) as $C =
F + G$,
where $F$ and $G$ are classes with $F^2 - 2p_a(F) = G^2 - 2p_a(G)
= 0$.
\item There are exactly $27$ relevant classes $[C]$ on $X_6$
with $C^2 - 2p_a(C) = 2$ and $p_a(C) = 1$.
Each such class $C$ is obtained from a unique relevant class $E$
with $E^2 - 2p_a(E) = -1$ by adding $E$ to the anticanonical class:
$C = -K + E$.
\end{enumerate}
\end{lemma}

The proof of the above lemma is left to the reader.
All of the statements can be easily shown by careful examination
of the table of relevant classes on $X_6$;
many of the statements are also elementary consequences of
intersection theory on rational surfaces.

\section{A formula for the quantum product}
\label{sectionformula}

It is obvious that the computation of the quantum product
depends on knowing the intersection numbers
\[
[A] \cdot (\alpha \otimes \beta \otimes \gamma) \in {\Bbb Z}
\]
for the relevant three-point loci $A$,
and for generators $\alpha$, $\beta$, and $\gamma$ of $H^*(\surf)$.

{}From the computations made in Examples \ref{phi0}-\ref{phiC2=2},
the degree of the closure of the locus ${\cal R}_{[C]}$ is an
important
number for the quantum product. We will denote this degree by
$d_{[C]}$:
\[
d_{[C]} = \text{degree of }\overline{{\cal R}_{[C]}}
\text{ in the projective space }|C|.
\]
The following lemma is then immediate from the computations
made in Examples \ref{phi0}-\ref{phiC2=2}.

\begin{lemma}
Let $\surf$ be a Del Pezzo surface,
and suppose that $C$ gives a relevant class on $\surf$
(or $C = 0$).
We denote by $(C \cdot C)$ the self-intersection of $C$
and by $p_a(C)$ its arithmetic genus.
Let $d_{[C]}$ denote the degree of the closure of the locus ${\cal
R}_{[C]}$
in the projective space $|C|$.
Then:
\begin{itemize}
\item in case $C=0$, the three-point class $[A_0]$
is the class of the small diagonal
and has real codimension $8$ in $\surf^3$.
The classes
$[\surf]\otimes[\surf]\otimes[p]$ and
$[D_1]\otimes[D_2]\otimes[\surf]$
(for divisors $D_i$) generate the complementary space
$H^4(\surf^3)$
(together with the classes obtained from these by permutations).
We have:
\begin{itemize}
\item[] $[A_0] \cdot ([\surf]\otimes[\surf]\otimes[p]) = 1$ and
\item[] $[A_0] \cdot ([D_1]\otimes[D_2]\otimes[\surf]) = (D_1\cdot
D_2)$.
\end{itemize}

\item in case $(C\cdot C) -2p_a(C)= -1$,
the image of ${\cal S}_{[C]}$ has real codimension $6$ in
$\surf^3$.
The classes
$[\surf]\otimes[D]\otimes[p]$ and $[D_1]\otimes[D_2]\otimes[D_3]$
(for divisors $D$, $D_i$) generate the complementary space
$H^6(\surf^3)$
(together with the classes obtained from these by permutations).
We have:
\begin{itemize}
\item[] $[A_{[C]}] \cdot ([\surf]\otimes[D]\otimes[p]) = 0$ and
\item[] $[A_{[C]}] \cdot ([D_1]\otimes[D_2]\otimes[D_3]) =
d_{[C]} (C \cdot D_1)(C \cdot D_2)(C \cdot D_3)$.
\end{itemize}
In this case $[A_{[C]}] = d_{[C]} [C]\otimes[C]\otimes[C] \in
H^6(\surf^3)$.

\item in case $(C\cdot C) -2p_a(C) = 0$,
the image of ${\cal S}_{[C]}$ has real codimension $4$ in
$\surf^3$.
The classes
$[\surf]\otimes[p]\otimes[p]$ and $[D_1]\otimes[D_2]\otimes[p]$
(for divisors $D_i$) generate the complementary space
$H^8(\surf^3)$
(together with the classes obtained from these by permutations).
We have:
\begin{itemize}
\item[] $[A_{[C]}] \cdot ([\surf]\otimes[p]\otimes[p]) = 0$ and
\item[] $[A_{[C]}] \cdot ([D_1]\otimes[D_2]\otimes[p]) =
d_{[C]} (C \cdot D_1)(C \cdot D_2)$.
\end{itemize}
In this case $[A_{[C]}] =
d_{[C]} [C]\otimes[C]\otimes[\surf] +
d_{[C]} [C]\otimes[\surf]\otimes[C] +
d_{[C]} [\surf]\otimes[C]\otimes[C]
\in H^4(\surf^3)$.

\item in case $(C\cdot C) -2p_a(C) = 1$,
the image of ${\cal S}_{[C]}$ has real codimension $2$ in
$\surf^3$.
The classes
$[D]\otimes[p]\otimes[p]$
(for divisors $D$) generate the complementary space
$H^{10}(\surf^3)$
(together with the classes obtained from these by permutations).
We have:
\begin{itemize}
\item[] $[A_{[C]}] \cdot ([D]\otimes[p]\otimes[p]) = d_{[C]} (C
\cdot D)$.
\end{itemize}
In this case $[A_{[C]}] =
d_{[C]} [C]\otimes[\surf]\otimes[\surf] +
d_{[C]} [\surf]\otimes[C]\otimes[\surf] +
d_{[C]} [\surf]\otimes[\surf]\otimes[C]
\in H^2(\surf^3)$.

\item in case $(C\cdot C) = 2$,
the image of ${\cal S}_{[C]}$ is all of $\surf^3$
(and therefore has codimension zero).
The class
$[p]\otimes[p]\otimes[p]$
generates the complementary space $H^{12}(\surf^3)$.
We have:
\begin{itemize}
\item[] $[A_{[C]}] \cdot ([p]\otimes[p]\otimes[p]) = d_{[C]}$
\end{itemize}
and $[A_{[C]}] = d_{[C]} [\surf]\otimes[\surf]\otimes[\surf] \in
H^0(\surf^3)$.
\end{itemize}
\end{lemma}

This gives the following descriptions
of the classes $\phi_{[C]}(\alpha,\beta)$:

\begin{corollary}
Let $\surf$ be a Del Pezzo surface,
and suppose that $C$ gives a relevant class on $\surf$
(or $C = 0$).
We denote by $(C \cdot C)$ the self-intersection of $C$
and by $p_a(C)$ the arithmetic genus.
Let $d_{[C]}$ denote the degree of the closure of the locus ${\cal
R}_{[C]}$
in the projective space $|C|$.
Then:
\begin{itemize}
\item in case $C=0$,
\begin{itemize}
\item[] $\phi_0([\surf],[\surf]) = [\surf]$,
\item[] $\phi_0([\surf],[p]) = [p]$,
\item[] $\phi_0([\surf],[D]) = [D]$ for a divisor $D$, and
\item[] $\phi_0([D_1],[D_2]) = (D_1\cdot D_2)[p]$ for divisors
$D_i$.
\item[] The classes $\phi_{[0]}([D],[p]) = \phi_{[0]}([p],[p]) =
0$.
\end{itemize}
\item in case $(C\cdot C) -2p_a(C) = -1$,
\begin{itemize}
\item[] $\phi_{[C]}([D_1],[D_2]) = d_{[C]} (C \cdot D_1)(C \cdot
D_2) [C]$.
\item[] If $\alpha$ and $\beta$ are homogeneous elements of
$H^*(\surf)$,
$\phi_{[C]}(\alpha,\beta) = 0$
unless both $\alpha$ and $\beta$ lie in $H^2(\surf)$.
\end{itemize}
\item in case $(C\cdot C) -2p_a(C) = 0$,
\begin{itemize}
\item[] $\phi_{[C]}([D_1],[D_2]) = d_{[C]} (C \cdot D_1)(C \cdot
D_2)[\surf]$,
 and
\item[] $\phi_{[C]}([D],[p]) = d_{[C]} (C \cdot D) [C]$.
\item[] $\phi_{[C]}([p],[p]) = 0$ and
$\phi_{[C]}([\surf],\beta) = 0$ for all $\beta$.
\end{itemize}
\item in case $(C\cdot C) -2p_a(C) = 1$,
\begin{itemize}
\item[] $\phi_{[C]}([p],[p]) = d_{[C]} [C]$ and
\item[] $\phi_{[C]}([D],[p]) = d_{[C]} (C \cdot D) [\surf]$.
\item[] $\phi_{[C]}([D_1],[D_2]) = 0$ and
$\phi_{[C]}([\surf],\beta) = 0$ for all $\beta$.
\end{itemize}
\item in case $(C\cdot C) -2p_a(C) = 2$,
\begin{itemize}
\item[] $\phi_{[C]}([p],[p]) = d_{[C]} [\surf]$.
\item[] If $\alpha$ and $\beta$ are homogeneous elements of
$H^*(\surf)$,
$\phi_{[C]}(\alpha,\beta) = 0$
unless both $\alpha$ and $\beta$ lie in $H^4(\surf)$.
\end{itemize}
\end{itemize}
\end{corollary}

Finally we deduce the formulas for the quantum product.

\begin{proposition}
\label{quantumproductformulas}
Let $\surf$ be a Del Pezzo surface.
\begin{enumerate}
\item The class $[\surf]$ is an identity for the quantum product.
\item For two divisors $D_1$ and $D_2$,
\[
[D_1] \qprod [D_2] = (D_1 \cdot D_2)[p] q^{[0]}
+ \sum\begin{Sb} E \text{ relevant}\\{E^2-2p_a(E) = -1}\end{Sb}
d_{[E]}(E \cdot D_1)(E \cdot D_2) [E] q^{[E]}
\]\[
+ \sum\begin{Sb} F \text{ relevant}\\{F^2-2p_a(F) = 0}\end{Sb}
d_{[F]} (F \cdot D_1)(F \cdot D_2) [\surf] q^{[F]}
\]
where the sum is taken over the linear systems (not over the curves
actually).
\item For a divisor $D$,
\[
[D] \qprod [p] =
\sum\begin{Sb} F \text{ relevant}\\{F^2-2p_a(F) = 0}\end{Sb}
d_{[F]} (F \cdot D) [F] q^{[F]}
+ \sum\begin{Sb} L \text{ relevant}\\{L^2-2p_a(L) = 1}\end{Sb}
d_{[L]} (L \cdot D) [\surf] q^{[L]}
\]
where the sum is again taken over the linear systems.
\item
\[
[p] \qprod [p] = \sum\begin{Sb} L \text{ relevant}\\{L^2-2p_a(L) =
1}\end{Sb}
  d_{[L]} [L] q^{[L]}
+ \sum\begin{Sb} C \text{ relevant}\\{C^2-2p_a(C) = 2}\end{Sb}
 d_{[C]} [\surf] q^{[C]}
\]
where the sum is again taken over the linear systems.
\end{enumerate}
\end{proposition}

\section{The relationship with ordinary cohomology}
\label{sectionordinary}

The effective cone $\Eff(\surf)$ in $H^2(\surf)$ is a proper cone,
in the sense that it contains no subgroups of $H^2(\surf)$.
Hence there is an ``augmentation'' ring homomorphism
\[
G: H^*_Q(\surf) \to H^*(\surf)
\]
defined by sending a quantum cohomology class $\sum_{[D]} \alpha_D
q^{[D]}$
to the coefficient $\alpha_0$ of the $q^{[0]}$ term.

\begin{proposition}
Let $\surf$ be a Del Pezzo surface.
Then the map $G$ is a ring homomorphism
from the quantum cohomology ring $H^*_Q(\surf)$ (with the quantum
product)
to the integral cohomology ring $H^*(\surf)$ (with the cup
product).
\end{proposition}

This is clear from the formulas for the quantum product given
in Proposition \ref{quantumproductformulas}.

\section{Examples}
\label{sectionexamples}

\begin{example}
Let $\surf = {\Bbb P}^2$, the complex projective plane.
\end{example}

Then $H^2(\surf) = {\Bbb Z}[L]$, where $[L]$ is the class of a
line.
The quantum products determining the multiplication are
\[
[L] \qprod [L] = [p]q^{[0]},
\]
\[
[L] \qprod [p] = [\surf]q^{[L]}, \text{ and }
\]
\[
[p] \qprod [p] = [L]q^{[L]}.
\]
We may identify $q^{[0]}$ and $[\surf]$ with $1$ and $q^{[L]}$ with
$q$;
if we write $\ell$ for the class $[L]$, the above relations are
that
$\ell^2 = [p]$, $\ell^3 = q$, and $\ell^4 = \ell q$.
Hence the quantum cohomology ring is isomorphic to
\[
H^*_Q({\Bbb P}^2) = {\Bbb Z}[\ell,q]/(\ell^3 - q).
\]

\begin{example}
Let $\surf = {\Bbb P}^1 \times {\Bbb P}^1$, the smooth quadric
surface.
\end{example}

Then $H^2(\surf) = {\Bbb Z}[F_1] \oplus {\Bbb Z}[F_2]$,
where the classes $[F_i]$ are those of the two rulings on $\surf$.
The quantum products determining the multiplication are
\[
[F_1] \qprod [F_1] = [\surf]q^{[F_2]},
\]
\[
[F_1] \qprod [F_2] = [p]q^{[0]},
\]
\[
[F_2] \qprod [F_2] = [\surf]q^{[F_1]},
\]
\[
[F_1] \qprod [p] = [F_2]q^{[F_2]},
\]
\[
[F_2] \qprod [p] = [F_1]q^{[F_1]}, \text{ and }
\]
\[
[p] \qprod [p] = [\surf]q^{[F_1+F_2]}.
\]
Denote $[F_i]$ by $f_i$, and $q^{[F_i]}$ by $q_i$.
These relations then become
$f_1^2 = q_2$, $f_1f_2 = [p]$, $f_2^2 = q_1$,
$f_1[p] = f_2q^2$, $f_2[p] = f_1q_1$, and $[p]^2 = q_1q_2$.
Hence the quantum cohomology ring is isomorphic to
\[
H^*_Q({\Bbb P}^1 \times {\Bbb P}^1) =
{\Bbb Z}[f_1,f_2,q_1,q_2]/(f_1^2-q_2,f_2^2-q_1).
\]

\begin{example}
Let $\surf = {\Bbb F}_1$, the blowup of the plane at one point.
\end{example}

Then $H^2(\surf) = {\Bbb Z}[E] \oplus {\Bbb Z}[F]$,
where $[E]$ is the class of the exceptional curve
and $[F]$ is the class of the fiber.
The only other class with a smooth rational curve
of self-intersection at most $2$ is the class $[L] = [E]+[F]$;
it has self-intersection $1$.
The quantum products determining the multiplication are
\[
[E] \qprod [E] = -[p]q^{[0]} + [E]q^{[E]} + [\surf]q^{[F]},
\]
\[
[E] \qprod [F] = [p]q^{[0]}  - [E]q^{[E]},
\]
\[
[F] \qprod [F] = [E]q^{[E]},
\]
\[
[E] \qprod [p] = [F]q^{[F]},
\]
\[
[F] \qprod [p] = [\surf]q^{[E]+[F]}, \text{ and }
\]
\[
[p] \qprod [p] = [L]q^{[E]+[F]}.
\]
Denote $[E]$ by $e$, $q^{[E]}$ by $q$,
$[F]$ by $f$, $q^{[F]}$ by $r$, and $[p]$ by $p$.
These relations then become
$e^2 = -p+eq+r$, $ef = p-eq$, $f^2 = eq$,
$ep = fr$, $fp = qr$, and $p^2 = (e+f)qr$.
We may eliminate $p$ from the generators since $p = ef+eq$;
after doing so, the first and third relations become
$e^2 = r - ef$ and $f^2 = eq$,
and the other relations formally follow from these two.
Hence the quantum cohomology ring is isomorphic to
\[
H^*_Q({\Bbb F}_1) = {\Bbb Z}[e,f,q,r]/
(e^2+ef-r,f^2-eq).
\]

\section{A functoriality property}
\label{sectionfunctor}

Let $\pi:\surf \to \surff$ be a general blowup
of a Del Pezzo surface $\surff$ at a
single point $p$,
with an exceptional curve $E$.
Then we have $\pi\times\pi\times\pi: \surf^3 \to \surff^3$.
Suppose that $C$ is an irreducible curve in $\surff$,
such that its cohomology class $[C]$
is relevant for the quantum cohomology of $\surff$.
We may assume that $p$ is not on $C$.
Then $[\pi^{-1}(C)] = \pi^*[C]$ as a class on $\surf$,
and has the same self-intersection and arithmetic genus as does
$[C]$;
therefore it is a relevant class for the quantum cohomology for
$\surf$.

In other words, the three-point classes $[A_{[C]}]$ on $\surff^3$
and $[A_{\pi^*[C]}]$ on $\surf^3$ are both defined.

\begin{lemma}
With the above notation and conditions,
if $[C] \neq 0$, then
\[
[A_{\pi^*[C]}] = {(\pi\times\pi\times\pi)}^*([A_{[C]}]).
\]
\end{lemma}

\begin{pf}
What is clear is that the nodal rational curve loci
${\cal R}_{[C]}$ and ${\cal R}_{\pi^*[C]}$ are birational,
so that $\pi\times\pi\times\pi$ induces a birational map
\[
{\cal S}_{\pi^*[C]} \to {\cal S}_{[C]}.
\]
since we are blowing up a point which is not on the general member
of ${\cal R}_{[C]}$.
This implies immediately that
\[
\pi\times\pi\times\pi([A_{\pi^*[C]}]) = [A_{[C]}]
\]
as classes on $\surff^3$.
We want to investigate the pull-back, not the image;
however this at least says that
$[A_{\pi^*[C]}]$ and ${(\pi\times\pi\times\pi)}^*([A_{[C]}])$
will differ only on the exceptional part of the map
$\pi\times\pi\times\pi$,
i.e.,  only over the fundamental locus
${\cal E} = (\{p\}\times \surff \times \surff) \cup
(\surff\times \{p\} \times \surff) \cup
(\surff\times \surff \times \{p\})$.

We can take the different cases up one by one.
If $(C \cdot C) - 2p_a(C) = -1$,
then the rational curve locus ${\cal R}_{[C]}$ is a finite set,
and by assumption the point $p$ is not on any member;
hence the loci in question are disjoint from the fundamental loci
for $\pi\times\pi\times\pi$, and there is nothing to prove.

Similarly if $(C\cdot C) -2p_a(C)= 2$, then $[A_{[C]}] =
d_{[C]}[\surff^3]$
and $[A_{\pi^*[C]}] = d_{\pi^*[C]}[\surf^3]$;
since $d_{[C]} = d_{\pi^*[C]}$, the result follows in this case.

Suppose that $(C \cdot C) -2p_a(C)= 0$, so that ${\cal R}_{[C]}$ is
a curve.
Let $C_1, \dots C_d$ be the members of ${\cal R}_{[C]}$ through $p$
(here $d = d_{[C]}$).
Then $[A_{[C]}]$ intersects the fundamental locus ${\cal E}$
exactly in the union of the loci
$(\{p\} \times C_j \times C_j) \cup
(C_j\times \{p\} \times C_j) \cup
(C_j\times C_j \times \{p\})$;
over this locus in $\surf^3$ is the union of the loci
$(E \times (E + \overline{C_j}) \times (E + \overline{C_j})) \cup
((E + \overline{C_j})\times \{p\} \times (E + \overline{C_j})) \cup
((E + \overline{C_j})\times (E + \overline{C_j}) \times \{p\})$
(where $(E + \overline{C_j})$ means the union of the exceptional
curve $E$
with the proper transform of $C_j$).
Note that this is a union of three-folds in $\surf^3$;
since in this case $[A_{\pi^*[C]}]$ is the class of a four-fold,
we have no extra contribution to the pull-back class.

The final case of $(C \cdot C) -2p_a(C) = 1$ is similar;
here the locus $[A_{[C]}]$ intersects the fundamental locus ${\cal
E}$
in a three-fold, over which lies a four-fold in $\surf^3$;
since $[A_{\pi^*[C]}]$ is the class of a five-fold,
again we have no extra contribution to the pull-back class.
\end{pf}

This implies the following.

\begin{corollary}
\label{pi*phi}
With the above notations, if $[C] \neq 0$, then
for any homogeneous classes $\alpha$ and $\beta$ in $H^*(\surff)$,
we have
\[
\pi^*(\phi_{[C]}(\alpha,\beta)) =
\phi_{[\pi^*C]}(\pi^*\alpha,\pi^*\beta).
\]
\end{corollary}

\begin{pf}
Choose a class $\gamma$ of the correct dimension on $\surff$,
and compute
\begin{eqnarray*}
\phi_{[\pi^*C]}(\pi^*\alpha,\pi^*\beta) \cdot \pi^*\gamma
&=& [A_{\pi^*[C]}] \cdot \pi^*\alpha\otimes \pi^*\beta \otimes
\pi^*\gamma \\
&=& {(\pi\times\pi\times\pi)}^*[A_{[C]}] \cdot
{(\pi\times\pi\times\pi)}^*(\alpha\otimes\beta\otimes\gamma \\
&=& [A_{[C]}] \cdot (\alpha\otimes\beta\otimes\gamma \\
&=& \phi_{[C]}(\alpha,\beta) \cdot \gamma.
\end{eqnarray*}
Moreover if $\phi_{[C]}(\alpha,\beta)$ is a class in $H^2(\surff)$,
and $E$ is the exceptional curve for the map $\pi$, then
\begin{eqnarray*}
\phi_{[\pi^*C]}(\pi^*\alpha,\pi^*\beta) \cdot E
&=& [A_{\pi^*[C]}] \cdot \pi^*\alpha\otimes \pi^*\beta \otimes E \\
&=& {(\pi\times\pi\times\pi)}^*[A_{[C]}] \cdot
\pi^*\alpha\otimes \pi^*\beta \otimes E \\
&=& 0.
\end{eqnarray*}
Hence as far as intersections go,
the class $\phi_{[\pi^*C]}(\pi^*\alpha,\pi^*\beta)$
is behaving exactly like the class
$\pi^*(\phi_{[C]}(\alpha,\beta))$.
However this class is defined in terms of its intersection
behaviour,
and so the equality as claimed holds.
\end{pf}

Note that we in any case have the formula
\[
\pi^*(\phi_{0}(\alpha,\beta)) = \phi_{0}(\pi^*\alpha,\pi^*\beta)
\]
since $\phi_{0}(\pi^*\alpha,\pi^*\beta) = \pi^*\alpha \cup
\pi^*\beta$
(this is cup product on $\surf$)
which is in turn equal to $\pi^*(\alpha \cup \beta)$
since $\pi^*$ is a ring homomorphism on ordinary cohomology.
Since $\pi^*(0) = 0$, we view this as the ``$[C] = 0$'' case
of Corollary \ref{pi*phi}.

Putting this together we derive the following version of
functoriality
for the quantum ring:

\begin{corollary}
\label{functoriality}
Let $\alpha$ and $\beta$ be ordinary cohomology classes in
$H^*(\surff)$.
Then for any relevant class $[C]$ in $H^2(\surff)$,
the $q^{\pi^*[C]}$-term of $\pi^*\alpha \qprod \pi^*\beta$
is equal to $\pi^*$ of the $q^{[C]}$-term of $\alpha\qprod\beta$.
\end{corollary}

Another way of saying this is to define the quantum pullback
\[
\pi_Q^*:H^*_Q(\surff) \to H^*_Q(\surf)
\]
by setting
\[
\pi^*_Q(\sum_{[D]} c_{[D]} q^{[D]} ) = \sum_{[D]} \pi^*(c_{[D]})
q^{\pi^*[D]}.
\]
This is NOT in general a ring homomorphism.
But the above corollary says that
for classes $\alpha$ and $\beta$ in $H^*(\surff)$,
the two quantum cohomology classes
\[
\pi^*(\alpha) \qprod \pi^*(\beta) \text{ and }
\pi^*_Q(\alpha \qprod \beta)
\]
differ only in the $q^{[D]}$ terms for those classes $[D]$ on
$\surf$
which are NOT pullbacks from $\surff$.
That is, they agree on all the $q^{[\pi^*C]}$ terms,
for any effective classes $C$ on $\surff$.

\section{Associativity of the quantum product for strict Del Pezzo
surfaces}
\label{sectionSDPassoc}

In this section we will use
the formulas of Proposition \ref{quantumproductformulas}
to check the associative law for the quantum product
for general strict Del Pezzo surfaces.
These are the surfaces ${\Bbb P}^2$,
${\Bbb F}_0 = {\Bbb P}^1 \times {\Bbb P}^1$,
and $X_n$ (the $n$-fold general blowup of ${\Bbb P}^2$)
for $n \leq 6$.
The first reduction is to note that it suffices to prove the
associative law
for the general surface $X_6$, the six-fold blowup of the plane.
(This is the general cubic surface in ${\Bbb P}^3$.)
This is due to the functoriality property
stated in Corollary \ref{functoriality};
if associativity holds on a blowup $\surf$ of a surface $\surff$,
then in fact it must hold on $\surff$.

By the tables of relevant classes on these surfaces
given in Section \ref{strictrelevantsection},
we note that all relevant classes have arithmetic genus at most
one.
(This is no longer true if one blows up $7$ general points in the
plane;
the class of quartics double at one point and passing through $6$
others
is relevant on $X_7$, and has arithmetic genus $2$.)
Hence for the strict Del Pezzo surface case,
the following lemma suffices to give us all the relevant $d_{[C]}$
numbers.

\begin{lemma}
Suppose that $\surf$ is a general rational surface,
and $[C]$ is a relevant class on $\surf$.
Then
\[
d_{[C]} = \left\{\begin{array}{cl}
1  & \text{ if }\;\; p_a(C) = 0, \text{ and } \\
12 & \text{ if }\;\; p_a(C) = 1.
\end{array}\right.
\]
\end{lemma}

\begin{pf}
If $p_a(C) = 0$, then the locus ${\cal R}_{[C]}$
of irreducible rational curves in the linear system $|C|$
is an open subset of $|C|$
(it is the subset parametrizing the smooth curves of $|C|$).
Hence its closure is the entire linear system $|C|$,
and therefore has degree one.

If $p_a(C) = 1$, then the locus ${\cal R}_{[C]}$
of irreducible rational curves in the linear system $|C|$
is an open subset of the discriminant locus of $|C|$.
Its degree is the number of irreducible nodal rational curves
in a general pencil of curves in $|C|$.
Such a general pencil, after blowing up the base points,
will give a fibration of elliptic curves on a rational surface;
the degree of ${\cal R}_{[C]}$ is the number of singular fibers
of this fibration.
This is $12$, by standard Euler number considerations.
\end{pf}

Let us begin with checking associativity for a triple product of
the form
$[p] \qprod [p] \qprod [D]$ for a divisor $D$.
For notational convenience let us define
\[
s(C) = (C \cdot C) - 2 p_a(C)
\]
for an irreducible curve $C$ in a relevant class $[C]$.
We have
\begin{eqnarray*}
([p] \qprod [p]) \qprod [D] &=&
(\sum_{s(L)= 1}  d_{[L]}[L] q^{[L]} +
\sum_{s(C)= 2} d_{[C]}[\surf] q^{[C]}) \qprod [D] \\
&=& \sum_{s(L) = 1} d_{[L]}([L]\qprod[D]) q^{[L]} +
\sum_{s(C)= 2} d_{[C]}[D] q^{[C]} \\
&=& \sum_{s(C)= 2} d_{[C]}[D] q^{[C]} +
\sum_{s(L)= 1} d_{[L]}(L \cdot D)[p] q^{[L]} + \\
&& + \sum_{s(L)= 1}\sum_{s(E)= -1}
d_{[L]} d_{[E]} (E \cdot L)(E \cdot D) [E] q^{[E+L]}\\
&& +
\sum_{s(L) = 1}\sum_{s(F)= 0}
d_{[L]} d_{[F]} (F \cdot L)(F \cdot D) [\surf] q^{[F+L]})
\end{eqnarray*}
while
\begin{eqnarray*}
[p] \qprod ([p] \qprod [D]) &=&
[p]\qprod (\sum_{s(F)= 0} d_{[F]}(F \cdot D) [F] q^{[F]}
+ \sum_{s(L)= 1} d_{[L]} (L \cdot D) [\surf] q^{[L]}) \\
&=& \sum_{s(F)= 0} d_{[F]}(F \cdot D) ([p]\qprod [F]) q^{[F]}
+ \sum_{s(L)= 1} d_{[L]} (L \cdot D) ([p]\qprod [\surf]) q^{[L]} \\
&=& \sum_{s(L)= 1} d_{[L]} (L \cdot D) [p] q^{[L]} +\\
&&+\sum_{s(F)= 0}\sum_{s(G)= 0}
d_{[F]}d_{[G]}(F \cdot D)(G \cdot F)[G] q^{[F+G]} + \\
&& + \sum_{s(F)= 0}\sum_{s(L)= 1}
d_{[F]}d_{[L]}(F \cdot D)(L \cdot F)[\surf] q^{[F+L]}.
\end{eqnarray*}

Comparing terms in the above two expressions we see that
this particular triple product is associative if and only if
\begin{eqnarray}
\label{ppD}
\sum_{s(C)= 2} d_{[C]} [D] q^{[C]} &+ &
\sum_{s(L)= 1}\sum_{s(E)= -1}
d_{[L]} d_{[E]} (E \cdot L)(E \cdot D) [E] q^{[E+L]} \\
&=& \sum_{s(F)= 0}\sum_{s(G)= 0}
d_{[F]}d_{[G]}(F \cdot D)(G \cdot F)[G] q^{[F+G]}. \nonumber
\end{eqnarray}
Of course only relevant classes are included in the above sums.

By Lemma \ref{X6relevantlemma},
if $L$ and $E$ are relevant classes with $s(L) = 1$ and $s(E) =
-1$,
then $0 \leq (L\cdot E) \leq 2$; if $(L \cdot E) = 0$
then there is no contribution in (\ref{ppD}).
If $(L\cdot E) \neq 0$ then $L+E$ is a relevant class with $s(L+E)
= 2$;
if $L \neq -K$ then $p_a(L+E) = (L\cdot E) - 1$,
and if $L = -K$ then $(L \cdot E) = 1$ and $p_a(L+E) = 1$.

Similarly if $F$ and $G$ are relevant classes with $s(F) = s(G) =
0$,
then $0 \leq (F\cdot G) \leq 2$; if $(F \cdot G) = 0$
then there is no contribution in (\ref{ppD}).
If $(F\cdot G) \neq 0$ then $F+G$ is a relevant class with $s(F+G)
= 2$
and $p_a(F+G) = (F \cdot G) - 1$.

Therefore the only terms which can appear in the associativity
formula
(\ref{ppD}) are those $q^{[C]}$ terms for relevant classes $[C]$
with $s(C) = 2$ (and $p_a(C) \leq 1$).
In fact we have the following.

\begin{lemma}
The associativity of the triple product $p \qprod p \qprod D$
on $X_6$ is equivalent to the following two formulas:
\begin{itemize}
\item[(a)] For every relevant class $C$ with $s(C) = 2$ and $p_a(C)
= 0$,
and every divisor $D$,
\[
[D] +
\sum\begin{Sb} E \\ s(E) = -1 \\ (C \cdot E) = 0 \end{Sb}
(D \cdot E) [E] =
\sum\begin{Sb} (F,G) \\ s(F) = s(G) = 0 \\ F+G = C \end{Sb}
(D \cdot F) [G].
\]
\item[(b)] For every relevant class $C$ with $s(C) = 2$ and $p_a(C)
= 1$,
\[
6[D] + 6(K+C \cdot D) [K+C] +
\sum\begin{Sb} E \\ s(E) = -1 \\(C \cdot E) = 1 \end{Sb}
(D \cdot E) [E] =
\sum\begin{Sb} (F,G) \\ s(F) = s(G) = 0 \\ F+G = C\end{Sb}
(D \cdot F) [G].
\]
\end{itemize}
\end{lemma}

\begin{pf}
As noted above, we may decompose (\ref{ppD}) into the $q^{[C]}$
terms
fixing a relevant class $C$ with $s(C) = 2$.  The two cases of the
lemma
correspond to the two possibilities for $p_a(C)$.

If $p_a(C) = 0$, then the only pairs $(L,E)$ with $s(L) = 1$ and
$s(E) = -1$
having $L + E = C$ must have $L \neq -K$ (and therefore $d_{[L]} =
1$).
Moreover $(L \cdot E) = 1$ (else $p_a(C) = 1$ by Lemma
\ref{X6relevantlemma}).
Hence $(C \cdot E) = (L + E \cdot E) = 1 - 1 = 0$.
Conversely for any class $E$ with $(C \cdot E) = 0$,
the class $L = C-E$ has $s(L) = 1$ and occurs in the sum.
Therefore this $(L,E)$ sum with $L+E = C$
is a sum over those $E$'s with $(C \cdot E) = 0$.
In this case $d_{[E]} = (L \cdot E) = 1$ also, so these
contributions
can be ignored.

Similarly, if $p_a(C) = 0$, then if $C = F+G$ with $s(F) = s(G) =
0$,
then $d_{[F]} = d_{[G]} = (F\cdot G) = 1$.
This then produces the equation of part (a).

Suppose then that $p_a(C) = 1$.
Then $d_{[C]} = 12$, and in the $(L,E)$ sum, $L = -K$ is a
possibility.
The $E$ that pairs with $L = -K$ is of course $E = K+C$,
and has $(E \cdot L) = (E \cdot -K) = 1$.
This gives a term $12(K+C\cdot D)[K+C]$ to the $(L,E)$ sum.
If $L \neq -K$ and $L+E = C$,
then by Lemma \ref{X6relevantlemma} we have $(L\cdot E) = 2$,
or, equivalently, $(C \cdot E) = 1$.
Conversely, any class $E$ with $s(E) = -1$ and $(C \cdot E) = 1$
occurs, and is paired with the class $L = C-E$.
Therefore again this sum can be written as a sum over such $E$'s,
each $E$ giving the term $2(E\cdot D)[E]$
(since $d_{[L]} = d_{[E]} = 1$ and $(L\cdot E) = 2$).

Finally, in the $(F,G)$ sum,
for two such classes to sum to a $C$ with $p_a(C) = 1$,
we must have $(F \cdot G) = 2$ by Lemma \ref{X6relevantlemma};
since $d_{[F]}=d_{[G]}=1$, each such pair $(F,G)$ contributes a
term
of the form $2 (F\cdot D)[G]$.

Dividing all terms by two produces the equation of part (b).
\end{pf}

It remains to prove these two formulas.
We begin with (a).

\begin{lemma}
Let $C$ be a relevant class on $X_6$ with $s(C) = 2$ and $p_a(C) =
0$.
Then for any divisor $D$ on $X_6$,
\[
[D] +
\sum\begin{Sb} E \\ s(E) = -1 \\ (C \cdot E) = 0 \end{Sb}
(D \cdot E) [E] =
\sum\begin{Sb} (F,G) \\ s(F) = s(G) = 0 \\ F+G = C \end{Sb}
(D \cdot F) [G].
\]
\end{lemma}

\begin{pf}
By Lemma \ref{X6relevantlemma},
the class $C$ can be written uniquely (up to order) as $C = F+G$,
with $s(G) = s(G) = 0$; if we do so, we see that
 the right-hand side of the above equation consists of only the two
terms
$(D \cdot F) [G] + (D \cdot G) [F]$.
Therefore we must actually show that
\[
[D] = (D \cdot F) [G] + (D \cdot G) [F] -
\sum\begin{Sb} E \\ s(E) = -1 \\ (C \cdot E) = 0 \end{Sb}
(D \cdot E) [E].
\]
The two pencils $|F|$ and $|G|$ on $X_6$ give a birational map
$\pi:X_6 \to {\Bbb F}_0$, realizing $X_6$ as a general five-fold
blowup
of ${\Bbb F}_0$.  The only curves $E$ on $X_6$ with $s(E) = -1$
which do not meet $C = F+G$ are the five exceptional curves for
this blowup;
call these five curves $E_1,\ldots,E_5$.
Note that the seven classes $[F],[G],[E_1],\ldots,[E_5]$
generate the Picard group over ${\Bbb Z}$;
the intersection matrix is unimodular.
Now the above formula is exactly the writing of the class $[D]$
in terms of these generators.
\end{pf}

Finally we address the equation (b).

\begin{lemma}
Let $C$ be a relevant class on $X_6$ with $s(C) = 2$ and $p_a(C) =
1$.
Then for any divisor $D$ on $X_6$,
\[
6[D] + 6(K+C \cdot D) [K+C] +
\sum\begin{Sb} E \\ s(E) = -1 \\(C \cdot E) = 1 \end{Sb}
(D \cdot E) [E] =
\sum\begin{Sb} (F,G) \\ s(F) = s(G) = 0 \\ F+G = C\end{Sb}
(D \cdot F) [G].
\]
\end{lemma}

\begin{pf}
Let $\hat{E}$ denote the class $K+C$;
we have $s(\hat{E}) = -1$ as noted above.
Any class $E$ with $s(E) = -1$ and $(C \cdot E) = 1$
must therefore have $(E \cdot \hat{E}) = 0$ and conversely;
therefore the $E$ sum above is a sum over those $E$'s
with $(E \cdot \hat{E}) = 0$.

On the other side, suppose that $F+G=C$ with $s(F)=s(G)=0$.
By Lemma \ref{X6relevantlemma}, we must have $(F \cdot G) = 2$,
and so $(F\cdot C) = (G \cdot C) = 2$.
Then $(F \cdot \hat{E}) = (F \cdot K+C) = (F \cdot K) + 2 = 0$
since $(F \cdot K) = -2$ by the genus formula.
Similarly $(G \cdot \hat{E}) = 0$.
Therefore in the two pencils $|F|$ and $|G|$,
the curve $\hat{E}$ occurs in a singular fiber of each.
Hence there are unique curves $E_F$ and $E_G$
with $s(E_F) = s(E_G) = -1$
such that $F = \hat{E} + E_F$ and $G = \hat{E} + E_G$.
Moreover $(E_F \cdot E_G) = 1$ since $(F \cdot G) = 2$.
Note that the three curves $\hat{E}$, $E_F$, and $E_G$ form a
triangle
on $X_6$ (considered as a cubic surface).

Conversely, given a triangle of curves $\hat{E}$, $E_F$, and $E_G$
with $s(E_F) = s(E_G)= -1$,
we obtain a unique pair $(F,G)$ with $s(F) = s(G) = 0$ and $F+G =
C$
by setting $G = \hat{E}+E_G$ and $F = \hat{E} + E_G$.
Therefore the $(F,G)$ sum above can be made into a sum over such
triangles;
for a fixed $\hat{E}$ there are five such \cite{beauville,reid}.
Therefore the equation in question may be written as
\[
6[D] + 6(\hat{E} \cdot D) [\hat{E}] +
\sum\begin{Sb} E \\ s(E) = -1 \\(\hat{E} \cdot E) = 0 \end{Sb}
(D \cdot E) [E] =
\sum\begin{Sb} \text{five triangles}\\ \hat{E}+E_F+E_G\end{Sb}
(D \cdot \hat{E}+E_F) [\hat{E}+E_G] + (D \cdot \hat{E}+E_G)
[\hat{E}+E_F].
\]

We first claim that the above equation holds when $D = \hat{E}$.
In this case the first two terms on the left side cancel,
while each term in the two sums are clearly zero.

Next we claim that the equation holds when $D = \hat{E} + E'$,
for any curve $E'$ with $s(E') = -1$ and $(\hat{E}\cdot E') = 1$.
In this case $(\hat{E} \cdot D) = 0$
so the second term of the equation drops out.
Of those $E$'s which satisfy $(\hat{E} \cdot E) = 0$,
there are exactly $8$ which meet $E'$
and contribute to the sum on the left-hand side;
these are exactly the other curves in the four triangles containing
$E'$
which do not involve $\hat{E}$.
These come in pairs, and if $E_1$ and $E_2$ form a pair,
then they contribute
$(D\cdot E_1)[E_1] + (D\cdot E_2)[E_2] = [E_1+E_2]$.
However each triangle is equivalent to $-K$,
so this pair's contribution may be written as $[-K - E']$;
hence this sum reduces to $-4[K+E']$.
Hence the entire left-hand side is equal to
$6[\hat{E} + E'] -4[K+E'] = 6[\hat{E}] - 4[K] + 2[E']$.

On the right-hand side,
if we have a triangle $E+E_F+E_G$,
$(D \cdot \hat{E}+E_F) = (E' \cdot \hat{E}+E_F) = 1 + (E'\cdot
E_F)$
and similarly for the $E_G$ term.
Now $E'$ is part of a triangle containing $\hat{E}$, say $\hat{E}
+ E' + E''$;
the other four triangles are disjoint from $E'$
(except for the curve $\hat{E}$).
For one of these four triangles, we obtain a contribution of
$[\hat{E}+E_G] + [\hat{E}+E_F] = [\hat{E} - K]$.
For the triangle with $E'$ and $E''$,
we have a contribution of
$2[\hat{E}+E']$.
Thus the right-hand sum reduces to $6[\hat{E}] - 4[K] + 2[E']$.

As noted above, this is equal to the left-hand side;
therefore the equation holds for this $D$.

The proof now finishes by remarking that the Picard group of $X_6$
is generated rationally by $\hat{E}$ and the classes $\hat{E}+E'$
considered above.
Since the equation is linear in $D$, and is true for these
generators,
it is true for all divisors $D$.
\end{pf}

Let us now address the associativity for a triple product of the
form
$[p]\qprod[D_1]\qprod[D_2]$ for divisors $D_i$.  We have
\[
[p]\qprod([D_1]\qprod[D_2]) =
\]
\[
{\renewcommand{\arraystretch}{1.5}
\begin{array}{ll}
= & [p]\qprod((D_1 \cdot D_2)[p] q^{[0]}
 + \sum_{s(E) = -1} d_{[E]}(E \cdot D_1)(E \cdot D_2) [E] q^{[E]}
 + \sum_{s(F) = 0} d_{[F]}(F \cdot D_1)(F \cdot D_2) [\surf]
q^{[F]}) \\
= & (D_1 \cdot D_2)[p]\qprod[p] q^{[0]}
+ \sum_{s(E) = -1} d_{[E]}(E \cdot D_1)(E \cdot D_2) [p]\qprod[E]
q^{[E]} \\
& + \sum_{s(F) = 0} d_{[F]}(F \cdot D_1)(F \cdot D_2)
[p]\qprod[\surf] q^{[F]}
\\
= & (D_1 \cdot D_2)(\sum_{s(L) = 1}d_{[L]}[L] q^{[L]} + \sum_{s(C)
= 2}
d_{[C]}[\surf] q^{[C]}) \\
& + \sum_{s(E) = -1} d_{[E]}(E \cdot D_1)(E \cdot D_2)
(\sum_{s(F) = 0} d_{[F]}(F \cdot E) [F] q^{[F]}
 +\sum_{s(L) = 1} d_{[L]}(L \cdot E) [\surf] q^{[L]}) q^{[E]} \\
& + \sum_{s(F) = 0}  d_{[F]}(F \cdot D_1)(F \cdot D_2) [p]q^{[F]}
\\
= & \sum_{s(L) = 1} d_{[L]}(D_1 \cdot D_2) [L] q^{[L]}
 + \sum_{s(C) = 2} d_{[C]}(D_1 \cdot D_2) [\surf] q^{[C]}  \\
& + \sum_{s(E) = -1}\sum_{F^2 = 0}
d_{[E]}d_{[F]} (E \cdot D_1)(E \cdot D_2)(F \cdot E) [F]
q^{[E]+[F]} \\
& + \sum_{s(E) = -1}\sum_{L^2 = 1}
d_{[E]}d_{[L]}(E \cdot D_1)(E \cdot D_2)(L \cdot E) [\surf]
q^{[E]+[L]} \\
& + \sum_{s(F) = 0} d_{[F]}(F \cdot D_1)(F \cdot D_2) [p]q^{[F]}
\end{array}
}
\]
while
\[
([p]\qprod[D_1])\qprod[D_2]) =
\]
\[
{\renewcommand{\arraystretch}{1.5}
\begin{array}{ll}
= & (\sum_{s(F) = 0} d_{[F]}(F \cdot D_1) [F] q^{[F]}
+ \sum_{s(L) = 1} d_{[L]}(L \cdot D_1) [\surf] q^{[L]}) \qprod
[D_2] \\
= & \sum_{s(F) = 0} d_{[F]}(F \cdot D_1) ([F]\qprod[D_2]) q^{[F]}
+ \sum_{s(L) = 1} d_{[L]}(L \cdot D_1) [D_2] q^{[L]} \\
= & \sum_{s(F) = 0} d_{[F]}(F \cdot D_1)
( (F \cdot D_2)[p] q^{[0]}
+ \sum_{s(E) = -1} d_{[E]}(E \cdot F)(E \cdot D_2) [E] q^{[E]} \\
& + \sum_{s(G) = 0} d_{[G]}(G \cdot F)(G \cdot D_2) [\surf] q^{[F]}
) q^{[F]} \\
& + \sum_{s(L) = 1} d_{[L]}(L \cdot D_1) [D_2] q^{[L]} \\
= & \sum_{s(F) = 0} d_{[F]}(F \cdot D_1) (F \cdot D_2)[p] q^{[F]}
\\
& + \sum_{s(F) = 0} \sum_{s(E) = -1}
d_{[F]}d_{[E]}(F \cdot D_1)(E \cdot F)(E \cdot D_2) [E] q^{[E+F]}
\\
& + \sum_{s(F) = 0} \sum_{s(G) = 0}
d_{[F]}d_{[G]}(F \cdot D_1)(G \cdot F)(G \cdot D_2) [\surf]
q^{[F+G]} \\
& + \sum_{s(L) = 1} d_{[L]}(L \cdot D_1) [D_2] q^{[L]} . \\
\end{array}
}
\]
Therefore associativity of this triple product is equivalent to the
identity
\[
{\renewcommand{\arraystretch}{1.5}
\begin{array}{ll}
& \sum_{s(L) = 1} d_{[L]}(D_1 \cdot D_2) [L] q^{[L]}
 + \sum_{s(C) = 2} d_{[C]}(D_1 \cdot D_2) [\surf] q^{[C]}  \\
& + \sum_{s(E) = -1}\sum_{s(F) = 0}
d_{[F]}d_{[E]}(E \cdot D_1)(E \cdot D_2)(F \cdot E) [F] q^{[E+F]}
\\
& + \sum_{s(E) = -1}\sum_{s(L) = 1}
d_{[L]}d_{[E]}(E \cdot D_1)(E \cdot D_2)(L \cdot E) [\surf]
q^{[E+L]} \\
=&\\
& \sum_{s(F) = 0} \sum_{s(E) = -1}
d_{[F]}d_{[E]}(F \cdot D_1)(E \cdot F)(E \cdot D_2) [E] q^{[E+F]}
\\
& + \sum_{s(F) = 0} \sum_{s(G) = 0}
d_{[F]}d_{[G]}(F \cdot D_1)(G \cdot F)(G \cdot D_2) [\surf]
q^{[F+G]} \\
& + \sum_{s(L) = 1} d_{[L]}(L \cdot D_1) [D_2] q^{[L]} .
\end{array}
}
\]
Comparing those terms with coefficients in $H^0(\surf)$ and those
in
$H^2(\surf)$, we see that
$ p \qprod (D_1 \qprod D_2) = (p \qprod D_1) \qprod D_2$ if and
only if the following two equations hold:
\begin{equation}
\label{pDD1}
{\renewcommand{\arraystretch}{1.5}
\begin{array}{ll}
& \sum_{s(L) = 1} d_{[L]}(D_1 \cdot D_2) [L] q^{[L]} \\
& + \sum_{s(E) = -1}\sum_{s(F) = 0}
d_{[F]}d_{[E]}(E \cdot D_1)(E \cdot D_2)(F \cdot E) [F] q^{[E+F]}
\\
=&\\
& \sum_{s(L) = 1}  d_{[L]}(L \cdot D_1) [D_2] q^{[L]} \\
& + \sum_{s(F) = 0} \sum_{s(E) = -1}
d_{[F]}d_{[E]}(F \cdot D_1)(E \cdot F)(E \cdot D_2) [E] q^{[E+F]}
\end{array}
}
\end{equation}
and
\begin{equation}
\label{pDD2}
{\renewcommand{\arraystretch}{1.5}
\begin{array}{ll}
& \sum_{s(C) = 2}  d_{[C]}(D_1 \cdot D_2) q^{[C]}  \\
& + \sum_{s(E) = -1}\sum_{s(L) = 1}
d_{[L]}d_{[E]}(E \cdot D_1)(E \cdot D_2)(L \cdot E) q^{[E+L]} \\
=&\\
& \sum_{s(F) = 0} \sum_{s(G) = 0}
d_{[F]}d_{[G]}(F \cdot D_1)(G \cdot F)(G \cdot D_2) q^{[F+G]} .\\
\end{array}
}
\end{equation}

\begin{lemma}
Equation (\ref{pDD2}) follows from (\ref{ppD}).
\end{lemma}

\begin{pf}
In fact, it is obtained from (\ref{ppD})
by setting $D = D_1$ and
dotting with $D_2$.
\end{pf}

The proof of associativity in the $p \qprod D_1 \qprod D_2$ case
now follows from the lemma below.

\begin{lemma}
Equation (\ref{pDD1}) holds for $X_6$.
\end{lemma}

\begin{pf}
We will show that (\ref{pDD1})
follows from two types of relations,
one of which holds generally for generic rational surfaces
and the other of which is special to the cubic surface.

Note that every relevant class $F$ with $s(F) = 0$
can be written uniquely as $F = -K - E_F$,
where $E_F$ is a relevant class with $s(E_F) = -1$
(see Lemma \ref{X6relevantlemma}).
Therefore for any relevant class $E$ with $s(E) = -1$,
we must have $0 \leq (F \cdot E) \leq 2$;
moreover $(F \cdot E) = 2$ if and only if $E = E_F$,
and $(F \cdot E) = 1$ if and only if $(E \cdot E_F) = 1$.
In this latter case $F+E$ is a relevant class with $s(F+E) = 1$
and $p_a(F+E) = 0$.

Therefore (\ref{pDD1}) can be analyzed by considering only these
types of
$q$-terms.  We begin by considering a term of the form $q^{[L]}$,
where $L$ is a relevant class with $s(L) = 1$ and $p_a(L) = 0$.

Note that in this case if $L=E+F$,
then $(E \cdot F) = 1$, $(E \cdot L) = 0$ and $(F\cdot L) = 1$.
Moreover all classes contributing to this term have $d=1$.
Thus considering the coefficent of $q^{[L]}$ in (\ref{pDD1})
gives the equation
\[
(D_1 \cdot D_2)[L] +
\sum\begin{Sb} s(E) = -1 \\ s(F) = 0 \\ E+F=L \end{Sb}
(E \cdot D_1)(E \cdot D_2)[F]
= (L \cdot D_1)[D_2] +
\sum\begin{Sb} s(E) = -1 \\ s(F) = 0 \\ E+F=L \end{Sb}
(E \cdot D_2)(F \cdot D_1)[E]
\]

Conversely we note that if $(L \cdot E) = 0$ for $E$ an exceptional
curve,
then by Riemann-Roch, $L\equiv E+F$ for some $F$.
Also if $E+F\equiv E'+F'$, then
$E\equiv E'$ and $F\equiv F'$ or $E$ and $E'$ are disjoint;
this follows from the fact that
$0 = (E' \cdot L) = (E' \cdot E) + (E' \cdot F) \geq (E' \cdot E)$
since $F$ moves in a pencil.
For each such class $L$,
it is easy to see that there are exactly $6$ classes $E$
with $s(E) = -1$ and $(E \cdot L) = 0$.
Therefore there is a disjoint basis of $\Pic(X)$,
$[L], [E_1], \ldots, [E_6]$ where $E_1, \dots , E_6$ are
exceptional curves,
and since $(L \cdot E_i) = 0$,
$L$ can be written as $F_i + E_i$ for each $i$.
In this basis, the equation above is equivalent to:

\[
(D_1 \cdot D_2)[L] - (L \cdot D_1)[D_2] =
\sum_i ((E_i \cdot D_2)(F_i \cdot D_1)[E_i] -
(E_i \cdot D_1)(E_i\cdot D_2)[F_i])
\]
and the right-hand side of this equation is equal to
\begin{eqnarray*}
&=& \sum_i ((E_i \cdot D_2)(F_i \cdot D_1)[E_i] -
(E_i \cdot D_1)(E_i\cdot D_2)[F_i]) \\
&=& \sum_i (E_i \cdot D_2)
((L \cdot D_1)-(E_i \cdot D_1)) [E_i] - (E_i \cdot D_1) [L-E_i]) \\
&=& \sum_i (E_i \cdot D_2)
(- (E_i \cdot D_1)[L] +
((L \cdot D_1)-(E_i \cdot D_1) + (E_i \cdot D_1)) [E_i]) \\
&=& (L \cdot D_1) \sum_i (E_i \cdot D_2)[E_i]
     - \sum_i (E_i \cdot D_2)(E_i \cdot D_1)[L]
\end{eqnarray*}
and so we must show that
\[
(D_1 \cdot D_2)[L] - (L \cdot D_1)[D_2] =
(L \cdot D_1) \sum_i (E_i \cdot D_2)[E_i]
     - \sum_i (E_i \cdot D_2)(E_i \cdot D_1)[L].
\]

On the other hand
$[D_2] = (L \cdot D_2)[L] - \sum_i(E_i \cdot D_2)[E_i]$.
Plugging this expression into the above equation,
we obtain an expression all in terms of the basis
$[L],[E_1],\ldots,[E_6]$.
The coefficients of $[E_i]$ on the two sides are obviously equal,
to $(L \cdot D_1)(E_i \cdot D_2)$.
Hence we must only check the coefficient of $[L]$, and so
the above equation follows from the identity
\[
(D_1 \cdot D_2) =
(L \cdot D_1)(L \cdot D_2) - \sum_i(E_i \cdot D_1)(E_i \cdot D_2)
\]
which is immediate from writing $[D_1]$ and $[D_2]$
in terms of the basis $[L],[E_1],\ldots,[E_6]$.

To complete the proof of the lemma,
we need to consider the coefficient of $q^{[-K]}$
in (\ref{pDD1});
$L = -K$ is the unique relevant class with $s=1$ and $p_a = 1$.
In the $E,F$ sums, we may sum over the $E$'s only, setting $F =
-K-E$;
noting that in this case $(E \cdot F) = 2$,
equating the coefficients of $q^{[-K]}$ in (\ref{pDD1})
and dividing by two gives

\begin{eqnarray}
6 (K \cdot D_1)[D_2] - 6(D_1 \cdot D_2)[K] &= &
\sum_{s(E)= -1}
((-K-E)\cdot D_1)(E \cdot D_2)[E] - (E \cdot D_2)(E \cdot
D_1)[-K-E]\nonumber
\\
&=& \sum_{s(E)= -1}
(-K \cdot D_1)(E \cdot D_2)[E] - (E \cdot D_2)(E \cdot D_1)[-K].
\label{pDD-K}
\end{eqnarray}

As we saw in the proof above for the associativity of the triple
product
$p \qprod p \qprod D$, $-K$ is
linearly equivalent to any triangle of exceptional curves;
moreover precisely ten exceptional curves meet any given
exceptional curve.
Thus for all exceptional curves $E'$,
\begin{eqnarray*}
5[-K] &=& 5[E'] + \sum_{(E \cdot E') = 1} [E] \\
&=& 6[E'] + \sum_E (E \cdot E')[E].
\end{eqnarray*}

Thus for all exceptional curves $E'$ and $E''$,
\[
6[E'] + \sum_E (E \cdot E')[E] = 6[E''] + \sum_E (E \cdot E'')[E]
\]
and so
\[
6[E'] = 6[E''] + \sum_E (E \cdot E'')[E] - \sum_E (E \cdot E')[E].
\]

Intersecting with $D_2$ and noting that $(-K \cdot E') = 1$ for all
$E'$,
we have
\[
6(E' \cdot D_2) = ( 6(E'' \cdot D_2) + \sum_E (E \cdot D_2)(E \cdot
E''))
(-K \cdot E') - \sum_E (E \cdot D_2)(E \cdot E').
\]

Here this equation holds for all exceptional curves $E'$,
which generate $\Pic(X_6)$, so
\[
6[D_2] = (6(E'' \cdot D_2) + \sum_E(E \cdot D_2)(E \cdot E''))[-K]
- \sum_E (E \cdot D_2)[E].
\]
Intersecting now with $D_1$ and again noting that $(-K \cdot
E'')=1$, we have
\[
(6(D_1 \cdot D_2)+ \sum_E (E \cdot D_2)(E \cdot D_1))(-K \cdot E'')
=
(-K \cdot D_1)(6(E'' \cdot D_2) + \sum_E(E \cdot D_2)(E \cdot
E'')).
\]
As this is true for all $E''$, which generate $\Pic(X_6)$, we see
that
\[
(6(D_1 \cdot D_2)+ \sum_E (E \cdot D_2)(E \cdot D_1))[-K]
= (-K \cdot D_1)(6[D_2] + \sum_E(E \cdot D_2)[E]).
\]
This can be re-written as
\begin{eqnarray*}
6(D_1 \cdot D_2)[-K] - 6(-K \cdot D_1)[D_2]
&=& \sum_E (-K \cdot D_1)(E \cdot D_2)[E] - (E \cdot D_1)(E \cdot
D_2)[-K]\\
&=& \sum_E (E \cdot D_2)((-K \cdot D_1)[E] - (E \cdot D_1)[-K]).
\end{eqnarray*}

This is exactly the desired equation (\ref{pDD-K}).
\end{pf}

To conclude our proof of associativity for the quantum product on
$X_6$,
we must deal with triple quantum products
of the form $D_1 \qprod D_2 \qprod D_3$ for divisors $D_i$.
Note that $d_{[C]} = 1$ whenever $[C]$ is a relevant class on $X_6$
with $s(c) \leq 0$;
we then compute:
\begin{eqnarray*}
D_1 \qprod (D_2 \qprod D_3) & =& D_1 \qprod ( (D_2 \cdot D_3)[p]q^0
+ \!\!\sum_{s(E)= -1} (E \cdot D_2)(E \cdot D_3)[E]q^E
+ \!\!\sum_{s(F)= 0} (F \cdot D_2)(F \cdot D_3)[X]q^F) \\
& = & (D_2 \cdot D_3)(D_1 \qprod [p])q^0
+ \sum_{s(E)= -1} (E \cdot D_2)(E \cdot D_3)(D_1 \qprod [E])q^E \\
&&+ \sum_{s(F)= 0} (F \cdot D_2)(F \cdot D_3)[D_1]q^F \\
& = & \sum_{s(F)= 0} (D_2 \cdot D_3)(F \cdot D_1)[F]q^F
+ \sum_{s(L)= 1} d_{[L]}(D_2 \cdot D_3)(L \cdot D_1)[X]q^L \\
&&+ \sum_{s(F)= 0} (F \cdot D_2)(F \cdot D_3)[D_1]q^F \\
&&+ \sum_{s(E)= -1} (E \cdot D_2)(E \cdot D_3)
( (E \cdot D_1)[p]q^E +
\sum_{s(E')= -1} (E' \cdot D_1)(E \cdot E')[E']q^{E+E'} \\
&&+ \sum_{s(F)= 0} (F \cdot D_1)(F \cdot E)[X]q^{E+F} ) \\
& = & \text{(dimension zero terms:)}
\sum_{s(E)= -1} (E \cdot D_1)(E \cdot D_2)(E \cdot D_3)[p]q^E \\
& + & \text{(dimension two terms:)}
\sum_{s(F)= 0} ( (D_2 \cdot D_3)(F \cdot D_1)[F] +
                 (F \cdot D_2)(F \cdot D_3)[D_1] ) q^F \\
&& + \sum_{s(E)=s(E')= -1}
     (E \cdot D_2)(E \cdot D_3)(E' \cdot D_1)(E \cdot
E')[E']q^{E+E'} \\
& + & \text{(dimension four terms:)}
\sum_{s(L)= 1} d_{[L]}(D_2 \cdot D_3)(L \cdot D_1)[X]q^L \\
&& + \sum\begin{Sb} s(E)= -1 \\ s(F)= 0 \end{Sb}
       (E \cdot D_2)(E \cdot D_3)(F \cdot D_1)(F \cdot E)[X]q^{E+F}
\\
\end{eqnarray*}

On the other hand
\begin{eqnarray*}
(D_1 \qprod D_2) \qprod D_3 & =& D_3 \qprod (D_1 \qprod D_2) \\
&=& \text{(dimension zero terms:)}
\sum_{s(E)= -1} (E \cdot D_3)(E \cdot D_1)(E \cdot D_2)[p]q^E \\
& + & \text{(dimension two terms:)}
\sum_{s(F)= 0} ( (D_1 \cdot D_2)(F \cdot D_3)[F] +
                 (F \cdot D_1)(F \cdot D_2)[D_3] ) q^F \\
&& + \sum_{s(E)=s(E')= -1}
     (E \cdot D_1)(E \cdot D_2)(E' \cdot D_3)(E \cdot
E')[E']q^{E+E'} \\
& + & \text{(dimension four terms:)}
\sum_{s(L)= 1} d_{[L]}(D_1 \cdot D_2)(L \cdot D_3)[X]q^L \\
&& + \sum\begin{Sb} s(E)= -1 \\ s(F)= 0 \end{Sb}
       (E \cdot D_1)(E \cdot D_2)(F \cdot D_3)(F \cdot E)[X]q^{E+F}
\\
\end{eqnarray*}
(This is obtained from the previous by permuting indices.)

Comparing terms, we see that the dimension zero terms are identical
and that the dimension four terms follow from Equation (\ref{pDD1})
(obtained by considering the $p \qprod D_1 \qprod D_2$ product)
intersected with $D_3$.

Comparing terms in dimension two, we need to show that
\begin{equation}
\label{DDD}
{\renewcommand{\arraystretch}{1.5}
\begin{array}{ll}
& \sum_{s(F)= 0}
   ( (D_1 \cdot D_2)(F \cdot D_3)[F] + (F \cdot D_1)(F \cdot
D_2)[D_3] )q^F \\
& + \sum_{s(E)=s(E')= -1}
  (E \cdot D_1)(E \cdot D_2)(E' \cdot D_3)(E' \cdot E) [E']
q^{E+E'} \\
=&\\
& \sum_{s(F)= 0}
   ( (D_2 \cdot D_3)(F \cdot D_1)[F] + (F \cdot D_2)(F \cdot
D_3)[D_1] ) q^F \\
& + \sum_{s(E)=s(E')= -1}
  (E' \cdot D_1)(E \cdot D_2)(E \cdot D_3)(E \cdot E')[E']q^{E+E'}.
\end{array}
}
\end{equation}

Now if $E$ and $E'$ are disjoint, then there is no contribution to
either side.  If $E = E'$, then the coefficients of $q^{E+E'}$ are
seen to be equal.  If $E$ meets $E'$, then $(E \cdot E')=1$ and
$E+E' \equiv F$ for some $F$.
Hence the equality above follows from the equality of the
coefficients of $q^F$
for a particular class $[F]$ with $s(F) = 0$.
Hence associativity is implied by the following lemma.

\begin{lemma}
For all relevant classes $[F]$ with $s(F)=0$,
\begin{eqnarray*}
&&(D_1 \cdot D_2)(F \cdot D_3)[F] + (F \cdot D_1)(F \cdot D_2)[D_3]
-(D_2 \cdot D_3)(F \cdot D_1)[F] - (F \cdot D_2)(F \cdot D_3)[D_1]
\\
&=&
\sum_{E+E' \equiv F}((E' \cdot D_1)(E \cdot D_2)(E \cdot D_3) - (E
\cdot D_1)(E \cdot D_2)(E' \cdot D_3))[E'].
\end{eqnarray*}
\end{lemma}

\begin{pf}
Our fixed curve $F$ gives $X$ a structure of ruled surface
with $F$ as general fiber.
Each $E+E'$ summing to $F$ is a reducible fiber
with respect to that ruling (there are five such reducible fibers).
Note that there is a basis, $[F],[G],[E_1],\dots,[E_5]$ of
$\Pic(X_6)$
such that $G$ has self-intersection  zero,
the $E_i$'s have self-intersection -1, $(F \cdot G) =1$,
and all other intersections are zero
(this is equivalent to $X_6$ having
${\Bbb F}_0 = {\Bbb P}^1\times {\Bbb P}^1$ as a minimal ruled
model).
We now write the right hand side of the above equation
in terms of the basis.
Note that $E'_i \equiv F - E_i$
and that the role of $E'$ is played by both $E'$ and $E$.
We have
\begin{eqnarray*}
&& \sum_{E+E' \equiv F}((E' \cdot D_1)(E \cdot D_2)(E \cdot D_3) -
(E
\cdot D_1)(E \cdot D_2)(E' \cdot D_3))[E'] \\
&=&
\sum_i ((F-E_i) \cdot D_1)(E_i \cdot D_2)(E_i \cdot D_3) - (E_i
\cdot D_1)(E_i \cdot D_2)((F-E_i) \cdot D_3))[F-E_i] \\
&& + \sum_i (E_i \cdot D_1)((F-E_i) \cdot D_2)((F-E_i) \cdot D_3)
-
((F-E_i) \cdot D_1)((F-E_i) \cdot D_2)(E_i \cdot D_3))[E_i].
\end{eqnarray*}

The coefficient of each $[E_i]$ in this expression
contains a sum of twelve products,
all of which cancel except for
$(E_i \cdot D_1)(F \cdot D_2)(F \cdot D_3) -
(F \cdot D_1)(F \cdot D_2)(E_i \cdot D_3)$.
The coefficient of the $[F]$ term is seen to be
$\sum_i  ( (F \cdot D_1)(E_i \cdot D_2)(E_i \cdot D_3)-
           (E_i \cdot D_1)(E_i \cdot D_2)(F \cdot D_3) )$.
Thus we may re-write the equation of the lemma as
\begin{eqnarray*}
&&(D_1 \cdot D_2)(F \cdot D_3)[F] + (F \cdot D_1)(F \cdot D_2)[D_3]
-(D_2 \cdot D_3)(F \cdot D_1)[F] - (F \cdot D_2)(F \cdot D_3)[D_1]
\\
&=& \\
&&
[F]\sum_i  ( (F \cdot D_1)(E_i \cdot D_2)(E_i \cdot D_3)-
           (E_i \cdot D_1)(E_i \cdot D_2)(F \cdot D_3) ) \\
&& + \sum_i ((E_i \cdot D_1)(F \cdot D_2)(F \cdot D_3) -
          (F \cdot D_1)(F \cdot D_2)(E_i \cdot D_3)) [E_i].
\end{eqnarray*}
To show that this linear equivalence is true,
it suffices to show equality when dotted with a basis of $\Pic(X)$.
It is easy to see
that equality holds when the expression above is dotted with $F$ or
any $E_i$.
Dotting with $G$, and recalling that $(F\cdot G) = 1$
and $(G \cdot E_i)= 0$ yields the following expression:
\begin{eqnarray*}
\label{mess}
&&(D_1 \cdot D_2)(F \cdot D_3) + (F \cdot D_1)(F \cdot D_2)(G\cdot
D_3)
-(D_2 \cdot D_3)(F \cdot D_1) - (F \cdot D_2)(F \cdot D_3)(G \cdot
D_1) \\
&=& \\
&&
\sum_i  ( (F \cdot D_1)(E_i \cdot D_2)(E_i \cdot D_3)-
           (E_i \cdot D_1)(E_i \cdot D_2)(F \cdot D_3) ).
\end{eqnarray*}
On the other hand, the divisors $D_j$, $j=1,2,3$
are written in terms of the basis as
\[
[D_j] \equiv (D_j \cdot G)[F] + (D_j \cdot F)[G] - \sum_i (D_j
\cdot E_i)[E_i].
\]
Substituting these expressions into $(D_1 \cdot D_2)$ and $(D_2
\cdot D_3)$
of (\ref{mess}) and collecting terms proves the result.
\end{pf}

This completes our analysis of the associativity of the quantum
product
for $X_6$, and therefore for all general strict Del Pezzo surfaces.
We have proved the following.
\begin{theorem}
The quantum product $\qprod$ is associative for
${\Bbb P}^2$, ${\Bbb F}_0$, and $X_1,\ldots,X_6$.
\end{theorem}

\section{Associativity in general}
\label{sectionassoc}

In this section we offer an algebro-geometric approach
to proving the associativity of the quantum product for a general
rational surface.
This approach avoids the reliance on perturbing
to a non-integrable almost complex structure
(see \cite{mcduff-salamon,ruan-tian});
we work with the existing complex/algebraic structure.

The associativity of the quantum product
is implied by checking associativity
for triple products of homogeneous generators for $H^*(\surf)$.
In other words, we must check that if $\alpha$, $\beta$, and
$\gamma$
are homogeneous classes in $H^*(\surf)$, then
\begin{equation}
\label{assoc1}
\alpha \qprod (\beta \qprod \gamma) =
(\alpha \qprod \beta) \qprod \gamma.
\end{equation}
We need not check the formula
when one of the constituents is the identity $[\surf]$,
or when they are all equal.

\begin{lemma}
\label{assoclemma}
The associativity of the quantum product
is equivalent to the following identity:
\[
\sum_{([C_1],[C_2]):[C_1+C_2]=[D]}
\phi_{[C_2]}(\alpha,\delta) \cdot \phi_{[C_1]}(\beta,\gamma)
=
\sum_{([C_1],[C_2]):[C_1+C_2]=[D]}
\phi_{[C_2]}(\gamma,\delta) \cdot \phi_{[C_1]}(\alpha,\beta).
\]
This identity must hold for all divisor classes $[D]$
and all homogeneous classes $\alpha$, $\beta$, $\gamma$, and
$\delta$
in $H^*(\surf)$.
\end{lemma}

\begin{pf}
Expanding the two sides of (\ref{assoc1}), we have
\begin{eqnarray*}
\alpha \qprod (\beta \qprod \gamma) &=&
\alpha \qprod ( \sum_{[C_1]} \phi_{[C_1]}(\beta,\gamma) q^{[C_1]})
\\
&=& \sum_{[C_1]} \sum_{[C_2]}
    \phi_{[C_2]}( \alpha,\phi_{[C_1]}(\beta,\gamma) ) q^{[C_1+C_2]}
\end{eqnarray*}
while
\begin{eqnarray*}
(\alpha \qprod \beta) \qprod \gamma &=&
( \sum_{[C_1]} \phi_{[C_1]}(\alpha,\beta) q^{[C_1]} )  \qprod
\gamma \\
&=& \sum_{[C_1]} \sum_{[C_2]}
    \phi_{[C_2]}(\phi_{[C_1]}(\alpha,\beta),\gamma ) q^{[C_1+C_2]}.
\end{eqnarray*}
For these to be equal,
they must have equal coefficients for all terms $q^{[D]}$.
Therefore associativity of the quantum product is equivalent to
having
\[
\sum_{([C_1],[C_2]):[C_1+C_2]=[D]}
\phi_{[C_2]}( \alpha,\phi_{[C_1]}(\beta,\gamma) )
=
\sum_{([C_1],[C_2]):[C_1+C_2]=[D]}
\phi_{[C_2]}(\phi_{[C_1]}(\alpha,\beta),\gamma )
\]
for all homogeneous classes $\alpha$, $\beta$, and $\gamma$ in
$H^*(\surf)$
and all divisor classes $[D]$.
The equality is equivalent to knowing that
for all homogeneous $\delta \in H^*(\surf)$,
the intersection products with $\delta$ are equal, i.e.,
\[
\sum_{([C_1],[C_2]):[C_1+C_2]=[D]}
\phi_{[C_2]}( \alpha,\phi_{[C_1]}(\beta,\gamma) ) \cdot \delta
=
\sum_{([C_1],[C_2]):[C_1+C_2]=[D]}
\phi_{[C_2]}(\phi_{[C_1]}(\alpha,\beta),\gamma ) \cdot \delta.
\]
(Here a dot product is taken to be zero
unless the codimensions of the classes are complementary.)

These intersection products can be computed on $\surf^3$,
and we then are requiring that
\[
\sum_{([C_1],[C_2]):[C_1+C_2]=[D]}
[A_{[C_2]}] \cdot
(\alpha\otimes\phi_{[C_1]}(\beta,\gamma)\otimes\delta)
=
\sum_{([C_1],[C_2]):[C_1+C_2]=[D]}
[A_{[C_2]}] \cdot (\phi_{[C_1]}(\alpha,\beta) \otimes\gamma \otimes
\delta).
\]

By the symmetry of the $[A]$-classes we may rewrite this as
\[
\sum_{([C_1],[C_2]):[C_1+C_2]=[D]}
[A_{[C_2]}] \cdot
(\alpha\otimes\delta\otimes\phi_{[C_1]}(\beta,\gamma))
=
\sum_{([C_1],[C_2]):[C_1+C_2]=[D]}
[A_{[C_2]}] \cdot ( \gamma \otimes
\delta\otimes\phi_{[C_1]}(\alpha,\beta)),
\]
which we may then reformulate as
\[
\sum_{([C_1],[C_2]):[C_1+C_2]=[D]}
\phi_{[C_2]}(\alpha,\delta) \cdot \phi_{[C_1]}(\beta,\gamma)
=
\sum_{([C_1],[C_2]):[C_1+C_2]=[D]}
\phi_{[C_2]}(\gamma,\delta) \cdot \phi_{[C_1]}(\alpha,\beta).
\]
This is the desired identity.
\end{pf}

Let $[C]$ be a relevant class on $\surf$,
and let ${\cal R}_{[C]}$ denote the locus of irreducible nodal
rational curves
in the linear system $|C|$.
Recall that $d_{[C]}$ is the degree of the closure $\overline{{\cal
R}_{[C]}}$.

Suppose that $[C]$ decomposes as $[C] = [C_1] +[C_2]$ where
$[C_1]$ and $[C_2]$ are relevant classes.
We may form the following locus
\begin{eqnarray*}
{\cal S}_{[C_1],[C_2]} &= \{&
(C_1,C_2,x_1,y_1,x_2,y_2,z)
\in {\cal R}_{[C_1]} \times {\cal R}_{[C_2]} \times \surf^5 \;|\;
\\
&& C_1 \text{ and } C_2 \text{ meet transversally at } z, \\
&& x_1 \text{ and } y_1 \text{ are smooth points of }C_1, \text{
and } \\
&& x_2 \text{ and } y_2 \text{ are smooth points of }C_2
\}.
\end{eqnarray*}

There is a natural map
\[
{\cal S}_{[C_1],[C_2]} \to \surf^6
\]
sending $(C_1,C_2,x_1,y_1,x_2,y_2,z)$ to $(x_1,y_1,z,x_2,y_2,z)$.
Call the image of the fundamental class $[A^6_{[C_1],[C_2]}]$.
Related to this is the map
\[
{\cal S}_{[C_1]} \times {\cal S}_{[C_2]} \to \surf^6
\]
sending a pair
$((C_1,x_1,y_1,z_1),(C_2,x_2,y_2,z_2))$ to
$(x_1,y_1,z_1,x_2,y_2,z_2)$.
The image of the fundamental class of this map is clearly
$[A_{[C_1]}] \otimes [A_{[C_2]}]$.
If we denote by $\pi_{36}:\surf^6 \to \surf^2$
the projection onto the third and sixth coordinates,
we see that
\[
[A^6_{[C_1],[C_2]}] = ([A_{[C_1]}] \otimes [A_{[C_2]}])
\cup \pi_{36}^*([\Delta])
\]
where $\Delta \subset \surf^2$ is the diagonal.

Finally consider the natural map
\[
{\cal S}_{[C_1],[C_2]} \to \surf^4
\]
sending $(C_1,C_2,x_1,y_1,x_2,y_2,z)$ to $(x_1,y_1,x_2,y_2)$.
Call the image of the fundamental class $[A^4_{[C_1],[C_2]}]$.

\begin{lemma}
With the above notation,
\[
\phi_{[C_1]}(\alpha,\beta) \cdot \phi_{[C_2]}(\gamma,\delta) =
[A^4_{[C_1],[C_2]}] \cdot
\alpha\otimes\beta\otimes\gamma\otimes\delta.
\]
\end{lemma}

\begin{pf}
Write $[\Delta] = \sum_i u_i \otimes v_i$ in $H^4(\surf^2)$,
where the $u_i$ and $v_i$ are classes in $H^*(\surf)$.
Denote by $\pi_{1245}:\surf^6 \to \surf^4$ the projection onto the
first, second, fourth, and fifth factors.
Then
\begin{eqnarray*}
\phi_{[C_1]}(\alpha,\beta) \cdot \phi_{[C_2]}(\gamma,\delta) &=&
[\phi_{[C_1]}(\alpha,\beta) \otimes \phi_{[C_2]}(\gamma,\delta)]
 \cdot [\Delta] \\
&=& \sum_i [\phi_{[C_1]}(\alpha,\beta) \otimes
\phi_{[C_2]}(\gamma,\delta)]
\cdot[u_i\otimes v_i] \\
&=& \sum_i (\phi_{[C_1]}(\alpha,\beta) \cdot u_i)
           (\phi_{[C_2]}(\gamma,\delta) \cdot v_i) \\
&=& \sum_i([A_{[C_1]}] \cdot \alpha\otimes\beta\otimes u_i)
    ([A_{[C_2]}] \cdot \gamma\otimes\delta\otimes v_i) \\
&=& \sum_i ([A_{[C_1]}] \otimes [A_{[C_2]}])\cdot
(\alpha\otimes\beta\otimes u_i \otimes \gamma\otimes\delta\otimes
v_i) \\
&=& ([A_{[C_1]}] \otimes [A_{[C_2]}])\cup
(\alpha\otimes\beta\otimes \surf \otimes \gamma\otimes\delta\otimes
\surf)
\cup \pi_{36}^*([\Delta]) \\
&=& ([A_{[C_1]}] \otimes [A_{[C_2]}])\cup
\pi_{1245}^*(\alpha\otimes\beta\otimes \gamma\otimes\delta)
\cup \pi_{36}^*([\Delta]) \\
&=& [A^6_{[C_1],[C_2]}] \cdot
\pi_{1245}^*(\alpha\otimes\beta\otimes \gamma\otimes\delta) \\
&=& [A^4_{[C_1],[C_2]}] \cdot
(\alpha\otimes\beta\otimes \gamma\otimes\delta)
\end{eqnarray*}
as claimed.
\end{pf}

Next we introduce the space
\begin{eqnarray*}
{\cal S}^4_{[C]} &= \{&
(C,x_1,x_2,x_3,x_4)
\in {\cal R}_{[C]} \times \surf^4 \;|\; \\
&& x_i \text{ are smooth points of }C
\}.
\end{eqnarray*}
There is a natural projection to $\surf^4$;
we denote the image of the fundamental class by $[A^4_{[C]}]$.

We note that we can consider the spaces ${\cal S}_{[C_1],[C_2]}$
as lying in the closure of the space ${\cal S}^4_{[C_1+C_2]}$,
via the natural map
\[
\psi:{\cal S}_{[C_1],[C_2]} \to \overline{ {\cal S}^4_{[C_1+C_2]}
}
\]
defined by sending $(C_1,C_2,x_1,y_1,x_2,y_2,z)$ to
$(C_1+C_2,x_1,y_1,x_2,y_2)$.
Actually, this map $\psi$ may be finite-to-one,
if $(C_1\cdot C_2) \geq 2$.
The image points represent the addition of an extra node to a curve
which already has arithmetic genus zero;
this then breaks the curve into two components.
One expects that the boundary of ${\cal S}^4_{[C]}$ in its closure
will have exactly the images of these loci
$\psi({\cal S}_{[C_1],[C_2]})$
with $[C_1]+[C_2] = [C]$
as components.

Now the space ${\cal S}^4_{[C]}$ has a cross-ratio function on it,
\[
\CR:{\cal S}^4_{[C]} \to {\Bbb P}^1,
\]
defined by sending $(C,x_1,x_2,x_3,x_4)$
to the cross-ratio
$(x_1-x_3)(x_2-x_4)/(x_1-x_4)(x_2-x_3)$.
(Here a coordinate is chosen on the normalization of $C$.)

Again one expects that this cross-ratio function will extend to the
closure
$\overline{ {\cal S}^4_{[C_1+C_2]} }$,
or at least to a model of the closure which is birational
on the boundary divisors $\psi({\cal S}_{[C_1],[C_2]})$
with $[C_1]+[C_2] = [C]$.

For distinct points, the cross-ratio takes values in
${\Bbb P}^1 - \{0,1,\infty\}$.
For the general coalescence of two of the four points,
the cross-ratio takes value $0$ when $x_1 = x_3$,
value $1$ when $x_1 = x_2$, and value $\infty$ when $x_1 = x_4$.

Therefore $\CR^{-1}(1)$ on the closure
should be the space where the first two points $x_1$ and $x_2$
come together; when this happens,
the curve will split, with $x_1$ and $x_2$ moving to points
on one curve and $x_3$ and $x_4$ lying on the other.
Therefore if we denote by $[A^4_{[C]}(\lambda)]$
the class of the image of $\CR^{-1}(\lambda)$ in $\surf^4$,
we have that
\[
[A^4_{[C]}(1)] = \sum\begin{Sb} ([C_1],[C_2]) \\ [C_1+C_2]=[C]
\end{Sb}
[A^4_{[C_1],[C_2]}]
\]
Therefore
\begin{eqnarray*}
\sum\begin{Sb} ([C_1],[C_2]) \\ [C_1+C_2]=[C] \end{Sb}
\phi_{[C_2]}(\alpha,\delta) \cdot \phi_{[C_1]}(\beta,\gamma)
&=& \sum\begin{Sb} ([C_1],[C_2]) \\ [C_1+C_2]=[C] \end{Sb}
[A^4_{[C_1],[C_2]}]
\cdot \beta\otimes\gamma\otimes\alpha\otimes\delta \\
&=& [A^4_{[C]}(1)] \cdot
\beta\otimes\gamma\otimes\alpha\otimes\delta.
\end{eqnarray*}

Now as $\lambda$ varies, the classes $[A^4_{[C]}(\lambda)]$
are rationally equivalent; hence the intersection product is the
same.
Therefore
\[
[A^4_{[C]}(1)] \cdot \beta\otimes\gamma\otimes\alpha\otimes\delta
=
[A^4_{[C]}(0)] \cdot \beta\otimes\gamma\otimes\alpha\otimes\delta.
\]

This class $[A^4_{[C]}(0)]$ can be analyzed
by studying the maps
\[
\tilde{\psi}:{\cal S}_{[C_1],[C_2]} \to \overline{ {\cal
S}^4_{[C_1+C_2]} }
\]
defined by sending $(C_1,C_2,x_1,y_1,x_2,y_2,z)$ to
$(C_1+C_2,x_1,x_2,y_1,y_2)$.
This is just the map $\psi$ above, followed
by a permutation of the four points;
but we see that if we denote the image of these classes by
$[\tilde{A}^4_{[C_1],[C_2]}]$, then
\[
[A^4_{[C]}(0)] = \sum\begin{Sb} ([C_1],[C_2]) \\ [C_1+C_2]=[C]
\end{Sb}
[\tilde{A}^4_{[C_1],[C_2]}]
\]
since the cross-ratio being zero represents when the first and
third
points come together, and this should be modelled by having them
split off
to one of the component curves (in this case $C_1$).

This class $[\tilde{A}^4_{[C_1],[C_2]}]$ is related to the original
version
$[A^4_{[C_1],[C_2]}]$ by the relation that
\[
[A^4_{[C_1],[C_2]}] \cdot
\beta\otimes\alpha\otimes\gamma\otimes\delta
=
[\tilde{A}^4_{[C_1],[C_2]}] \cdot
\beta\otimes\gamma\otimes\alpha\otimes\delta
\]
because of the permutation
which relates the maps $\psi$ and $\tilde{\psi}$.
Hence
\begin{eqnarray*}
[A^4_{[C]}(0)] \cdot \beta\otimes\gamma\otimes\alpha\otimes\delta
&=&
\sum\begin{Sb} ([C_1],[C_2]) \\ [C_1+C_2]=[C] \end{Sb}
[\tilde{A}^4_{[C_1],[C_2]}] \cdot
\beta\otimes\gamma\otimes\alpha\otimes\delta \\
&=& \sum\begin{Sb} ([C_1],[C_2]) \\ [C_1+C_2]=[C] \end{Sb}
[A^4_{[C_1],[C_2]}] \cdot
\beta\otimes\alpha\otimes\gamma\otimes\delta \\
&=& \sum\begin{Sb} ([C_1],[C_2]) \\ [C_1+C_2]=[C] \end{Sb}
\phi_{[C_1]}(\beta,\alpha) \cdot \phi_{[C_2]}(\gamma,\delta)
\end{eqnarray*}

We conclude that
\[
\sum\begin{Sb} ([C_1],[C_2]) \\ [C_1+C_2]=[C] \end{Sb}
\phi_{[C_2]}(\alpha,\delta) \cdot \phi_{[C_1]}(\beta,\gamma)
=
\sum\begin{Sb} ([C_1],[C_2]) \\ [C_1+C_2]=[C] \end{Sb}
\phi_{[C_1]}(\beta,\alpha) \cdot \phi_{[C_2]}(\gamma,\delta)
\]
which is equivalent to the associative law for the quantum product
by Lemma \ref{assoclemma}, using the symmetry of the
$\phi$-classes.

The reader will note that the approach given above
to the proof of associativity
relies on the existence of a model of the closure of
the spaces ${\cal S}^4_{[C]}$,
which has rather nice properties:
the boundary is well-understood in terms of splittings of $C$ as
$C = C_1 + C_2$,
and the cross-ratio function extends nicely to it.
The existence of such a model we only conjecture,
and have not attempted to construct it in this paper.

\section{Enumerative consequences of associativity}
\label{sectionenum}

We will now extract several enumerative consequences
from the associative law for the quantum product
which have been noted by Kontsevich and Manin in
\cite{kontsevich-manin}.

Let us change notation somewhat and introduce the integer
\[
k(C) = k([C]) =\begin{cases}
(-K\cdot C) & \text{ if the class $[C]$ is relevant on $\surf$}\\
0 & \text{if $[C]$ is not relevant.}
\end{cases}
\]
We note that by the adjunction formula, if $[C]$ is a relevant
class,
then $k(C) = s(C) + 2$
(recall that $s(C) = (C\cdot C) - 2p_a(C)$).
Moreover $k(-)$ is linear in relevant
classes:
\[
k(C_1+C_2) = k(C_1)+k(C_2).
\]

Also, the subscript notation for the degree $d_{[C]}$ of the
rational curve
locus ${\cal R}_{[C]}$ in the linear system $|C|$
is too cumbersome; we will switch notation and call this degree
$N(C)$
(following the notation in \cite{kontsevich-manin}).

The equations which were seen in Section \ref{sectionSDPassoc}
to be equivalent to the associative law
for the quantum product for strict Del Pezzos
are actually equivalent to associativity for any rational surface
$\surf$.
These were (\ref{ppD}), (\ref{pDD1}), and (\ref{DDD}).
With the above notation, they can be written as

\begin{eqnarray}
\label{assocc1}
\sum_{k(C)= 4} N(C) [D] q^{[C]} &+ &
\sum_{k(L)= 3}\sum_{k(E)= 1}
N(L) N(E) (E \cdot L)(E \cdot D) [E] q^{[E+L]} \\
&=& \sum_{k(F)= 2}\sum_{k(G)= 2}
N(F)N(G)(F \cdot D)(G \cdot F)[G] q^{[F+G]} \nonumber
\end{eqnarray}
for any divisor class $[D]$,
\begin{equation}
\label{assoc2}
{\renewcommand{\arraystretch}{1.5}
\begin{array}{ll}
& \sum_{k(L) = 3} N(L)(D_1 \cdot D_2) [L] q^{[L]} \\
& + \sum_{k(E) = 1}\sum_{k(F) = 2}
N(F)N(E)(E \cdot D_1)(E \cdot D_2)(F \cdot E) [F] q^{[E+F]} \\
=&\\
& \sum_{k(L) = 3}  N(L)(L \cdot D_1) [D_2] q^{[L]} \\
& + \sum_{k(F) = 2} \sum_{k(E) = 1}
N(F)N(E)(F \cdot D_1)(E \cdot F)(E \cdot D_2) [E] q^{[E+F]}
\end{array}
}
\end{equation}
for any divisor classes $[D_1]$ and $[D_2]$, and
\begin{equation}
\label{assoc3}
{\renewcommand{\arraystretch}{1.5}
\begin{array}{ll}
& \sum_{k(F)= 2} N(F)
   ( (D_1 \cdot D_2)(F \cdot D_3)[F] + (F \cdot D_1)(F \cdot
D_2)[D_3] )q^F \\
& + \sum_{k(E)=k(E')= 1} N(E)N(E')
  (E \cdot D_1)(E \cdot D_2)(E' \cdot D_3)(E' \cdot E) [E']
q^{E+E'} \\
=&\\
& \sum_{k(F)= 2}N(F)
   ( (D_2 \cdot D_3)(F \cdot D_1)[F] + (F \cdot D_2)(F \cdot
D_3)[D_1] ) q^F \\
& + \sum_{k(E)=k(E')= 1} N(E)N(E')
  (E' \cdot D_1)(E \cdot D_2)(E \cdot D_3)(E \cdot E')[E']q^{E+E'}.
\end{array}
}
\end{equation}
for any divisor classes $[D_1]$, $[D_2]$, and $[D_3]$.
(Note that (\ref{DDD}) had no $N(-)$ numbers; all were equal to one
in the case of $X_6$, but in general they must be included of
course.)

Suppose that $\surf = X_n$, the general $n$-fold blowup of the
plane.
We take as a basis for $\Pic(\surf)$ the classes
$[H]$, $[E_1], \dots, [E_n]$,
where $[H]$ is the class of the pullback of a line from ${\Bbb
P}^2$,
and $E_i$ is the exceptional curve over the $i$-th point $p_i$
which is blown up.
In this case every divisor class can be written as
\[
[D] = d[H] - \sum_{i=1}^n m_i [E_i]
\]
which we will abbreviate to $[D] = (d;m_1,\dots,m_n)$.
Note then that the anticanonical class $[-K] = (3;1^n)$
where we use the exponential notation for repeated $m_i$'s,
as is rather standard.
Hence if $[D] = (d;m_1,\dots,m_n)$ then $(-K\cdot D) = 3d-\sum_i
m_i$.

With this notation we see that quantum cohomology
is a completely numerically based theory.
A class $[D] = (d;m_1,\dots,m_n)$ is relevant if and only if
$d^2 + 3d \geq \sum_i m_i^2 + \sum_i m_i$
and $k(D) \leq 4$.
(The first condition is that the expected dimension of $|D|$ is
non-negative,
so that there will be curves in $|D|$;
the second is the relevance condition, that the dimension of the
locus
of rational curves in $|D|$ is not more than $3$.)

Suppose that we take the associativity condition (\ref{assocc1}),
set $D = H$, and dot with $H$: we obtain
\begin{eqnarray}
\label{assoc1D=HdotH}
\sum_{k(C)= 4} N(C) q^{[C]} &= &
\sum_{k(F)= 2}\sum_{k(G)= 2}
N(F)N(G)(F \cdot H)(G \cdot F)(G \cdot H) q^{[F+G]} \\
&& - \sum_{k(L)= 3}\sum_{k(E)= 1}
N(L) N(E) (E \cdot L){(E \cdot H)}^2q^{[E+L]} \nonumber
\end{eqnarray}
We want to use this to develop a recursive formula for the degrees
$N(C)$
if possible.

The $q^{[C]}$ terms of the (\ref{assoc1D=HdotH}) are
\begin{eqnarray}
\label{assoc1qC}
N(C) &= &
\sum\begin{Sb} (F,G) \\ k(F)=k(G)=2 \\ F+G \equiv C \end{Sb}
N(F)N(G)(F \cdot H)(G \cdot F)(G \cdot H) \\
&& - \sum\begin{Sb} (E,L) \\ k(E)=1, k(L) = 3 \\ E+L \equiv C
\end{Sb}
N(L) N(E) (E \cdot L){(E \cdot H)}^2 \nonumber
\end{eqnarray}

To relate this recursive formula to those of Kontsevich and Manin
(Claims 5.2.1 and 5.2.3b of \cite{kontsevich-manin}), we write this
as
\begin{eqnarray}
\label{assoc2qC}
 & & \\
N(C) &= &
\sum\begin{Sb} (C_1,C_2) \\ C_1+C_2 \equiv C \end{Sb}
N(C_1)N(C_2)(C_1 \cdot C_2)(H \cdot C_1)
( (H \cdot C_2)\delta_{k(C_1)-2} - (H \cdot C_1)\delta_{k(C_1)-1}
)  \nonumber
\end{eqnarray}
where $\delta_n = \begin{cases} 1 &  \text{ if $n=0$ } \\
0 & \text{ if $n\ne 0$ }  \end{cases}$.  This expression is
equivalent to Claim 5.2.3b of \cite{kontsevich-manin} in the case
for which $k(C) = 4$, assuming that their convention for
$\left( \begin{array}{c} 0 \\ n \end{array} \right)$ is that
$\left( \begin{array}{c} 0 \\ n \end{array} \right) = \delta_n$.
Note that \ref{assoc2qC} is not valid when $k(C) \ne 4$, as may be
seen when $C \equiv -K_X$ on $X=X_6$. It is unclear what 5.2.3b of
\cite{kontsevich-manin} means in this case, since  $(-K \cdot C) =
3$ and so 5.2.3b involves terms of the form $\left(
\begin{array}{c} -1 \\ n \end{array} \right)$ where $n \le 0$. Note
also
that the
same $C$ is indecomposible in the semi-group of numerically
effective curves and yet has $N(C)=12$, contrary to the expectation
expressed in 5.2.3b of \cite{kontsevich-manin}.

Our next goal is to compute $N(d)$, the degree of the locus of
rational curves of degree $d$ in the plane.
This is not a relevant class on the plane, unless $d = 1$.
Since forcing a curve to pass through a generically chosen point
is a linear condition on the linear system, we have
\[
N(d;m_1,\dots,m_n,1) = N(d;m_1,\dots,m_n).
\]
Hence by induction we have that
\[
N(d) = N(d;1^{3d-4})
\]
in particular.
Now on the surface $X_{3d-4}$,
the class $C = (d;1^{3d-4})$ is a relevant class;
in fact $k(d;1^{3d-4}) = 4$ and so \ref{assoc1qC} may be used to
compute $N(d)$. We now need to understand those $k(E) = 1$ and
$k(L) = 3$ classes which sum to $C$,
and those $k(F) = k(G) = 2$ classes which sum to $C$.

First consider the class $E_i$ itself, which has $k(E) = 1$.
However $(E_i \cdot H) = 0$,
so these $k(E)=1$ classes do not contribute to the recursive
formula of
(\ref{assoc1qC}).

Hence we may assume $E$ is a relevant class with all $m_i$'s
non-negative.
In this case since $C = (d;1^{3d-4})$, all $m_i$'s for $E$
(and for the complementary $k=3$ class $L$) must be $0$ or $1$.
Therefore $E = (e;1^{3e-1})$ for some $e$ with $1 \leq e \leq d-1$,
where this notation means that $3e-1$ of the $m_i$'s are $1$,
and all the others are zero.
(There are
{\small $\left(\begin{array}{c} 3d-4 \\ 3e-1 \end{array}\right)$}
such classes.)
The complementary class $L$ is of the form $L = (d-e;1^{3d-3e-3})$
where the $1$'s occur in the complementary positions.
Note that with this notation
$(E\cdot L) = e(d-e)$ and $(E \cdot H) = e$,
with $N(E) = N(e)$ and $N(L) = N(d-e)$;
so the second sum above reduces to
\[
\sum_{E+L\equiv C}
N(L) N(E) (E \cdot L){(E \cdot H)}^2
=
\sum_{e=1}^{d-1}
\left(\begin{array}{c} 3d-4 \\ 3e-1 \end{array}\right)
N(e) N(d-e) (d-e) e^3.
\]

Now suppose that $F$ is a $k=2$ class;
again all its multiplicity numbers $m_i$ must be zero or one,
and so $F$ must have the form $F = (e;1^{3e-2})$
for some $e$ with $1 \leq e \leq d-1$;
there are
{\small $\left(\begin{array}{c} 3d-4 \\ 3e-2 \end{array}\right)$}
such classes.
The complementary class $G$ is $G = (d-e;1^{3d-3e-2})$,
where again the $1$'s occur in the complementary positions.
Note that with this notation
$(F\cdot G) = e(d-e)$, $(F \cdot H) = e$, and $(G \cdot H) = d-e$;
also $N(F) = N(e)$ and $N(G) = N(d-e)$.
Hence the first sum above reduces to
\[
\sum_{F+G\equiv C}
N(F)N(G)(F \cdot H)(G \cdot F)(G \cdot H)
=
\sum_{e=1}^{d-1}
\left(\begin{array}{c} 3d-4 \\ 3e-2 \end{array}\right)
N(e) N(d-e) e^2{(d-e)}^2.
\]

Collecting terms gives the following recursion relation for the
degrees $N(d)$:
\begin{equation}
\label{km5.2.1}
N(d) = \sum_{e=1}^{d-1} e^2(d-e)
\left[ (d-e) \left(\begin{array}{c} 3d-4 \\ 3e-2 \end{array}\right)
- e \left(\begin{array}{c} 3d-4 \\ 3e-1 \end{array}\right) \right]
N(e) N(d-e).
\end{equation}
This is exactly the enumerative prediction made by Kontsevich and
Manin
(Claim 5.2.1 of \cite{kontsevich-manin}).

We note here that from our point of view
this prediction follows from the associativity of the quantum
product
for arbitrarily large blowups of the plane; it is not enough to
know
it just for the plane.

To further illustrate the geometric and enumerative significance of
associativity of quantum cohomology, we return to $X_6$.

\begin{proposition}
Associativity of quantum cohomology on strict Del Pezzo surfaces is
equivalent to the fact that there are 27 exceptional curves on
$X_6$, each of which meets precisely 10 others.
\end{proposition}
\begin{pf}
In Section \ref{sectionSDPassoc}, associativity for $X_6$ (and so
all strict Del Pezzos) was shown using the fact that the number and
mutual disposition of exceptional curves are as above.  To complete
our proof, it suffices to show that if $\hat E$ is an exceptional
curve on $X_6$, $m$ is the number of other exceptional curves
meeting $\hat E$ and $e$ is the total number of exceptional curves,
then $m=10$ and $e=27$.  Consider the associativity relation
\ref{assoc2} with $ L \equiv D_1 \equiv D_2 \equiv -K $. Allowing
ourselves the knowledge that $N(F)=N(E)=1$ for all cases for which
$F+E \equiv -K$ and collecting terms, \ref{assoc2} implies that
\begin{eqnarray}
0 & = & \sum\begin{Sb} (E,F) \\ k(F)=2, k(E)=1 \\ E+F \equiv -K
\end{Sb} (E \cdot -K)( (F \cdot -K)[E] - (E \cdot -K)[F] )
\nonumber \\
 & = & \sum_{k(E)=1} 2[E] - [-K-E] \nonumber \\
 & = & \sum_{k(E)=1} ( 3[E] + [K] ) \nonumber
\end{eqnarray}
Intersecting with our fixed exceptional curve $\hat E$, we have $
0 = \sum_E ( 3 (E \cdot \hat E) - 1)$ and so $0 = -4 + 2m - (e-
(m+1))$ and thus $e = 3m - 3$.

On the other hand, applying the associativity relation
\ref{assocc1}
with $C \equiv -K + \hat E$, we get for all divisors $D$
\[
12[D]  =  \sum_{F+G=C} 1 \cdot 1 \cdot (F \cdot D) \cdot 2 [G]
            - \sum\begin{Sb} L+E \equiv -K+\hat E \\ p_a(L) = 0
\end{Sb}
              1 \cdot 1 \cdot 2 \cdot (E \cdot D)[E]
            - 12 (\hat E \cdot D)[\hat E]
\]
which is equivalent, after noting that $F+G \equiv C$ if and only
if $F \equiv \hat E + E_F$ for some exceptional curve $E_F$, to
\[
6( [D] + (\hat E \cdot D)[\hat E] ) =
   \sum\begin{Sb} \hat E + E_F \\ \hat E + E_G \end{Sb}
    ( (\hat E + E_F) \cdot D)[\hat E + E_G]
   - \sum_{(E \cdot \hat E) = 0} (E \cdot D)[E].
\]
Letting $D \equiv -K$ and intersecting with $-K$ yields $24 = 4m -
(e-(m+1))$ and so $e = 5m -23$, which combined with our previous
equation yields $e=27$ and $m=10$.
\end{pf}

Note finally that
$N(d;m_1, \dots , m_n)$ is invariant under permutations and Cremona
transformations.  The former is obvious and the latter follows from
the fact that $k$ and the decomposition of a curve into sums of
curves are invariant under Cremona transformations.
Also of course the arithmetic genus of a class is
invariant under symmetries and Cremona transformations.
It is tempting to conjecture that the number $N(d;m_1, \dots ,
m_n)$
depends only on the genus.
However a recent computation of A. Grassi \cite{grassi}
shows that this is not the case in general for classes
with arithmetic genus at least $2$.

\end{document}